\documentclass[11pt]{amsart}
\usepackage{amscd,amssymb,graphics,color,hyperref}
\usepackage{mathrsfs}

\newtheorem{theorem}{Theorem}[section]
\newtheorem{lemma}[theorem]{Lemma}
\newtheorem{proposition}[theorem]{Proposition}

\newtheorem*{Thm.A}{Theorem A}
\newtheorem*{Thm.B}{Theorem B}
\newtheorem*{Thm.C}{Theorem C}

\theoremstyle{definition}
\newtheorem{definition}[theorem]{Definition}
\newtheorem{example}[theorem]{Example}

\theoremstyle{remark}
\newtheorem{remark}[theorem]{Remark}
\newtheorem{remarks}[theorem]{Remarks}

\newcommand{\NN} {\mathbb{N}}
\newcommand{\ZZ} {\mathbb{Z}}

\newcommand{\RR} {\mathbb{R}}
\newcommand{\CC} {\mathbb{C}}

\newcommand{\PP} {\mathbb{P}}

\newcommand {\shC}  {\mathcal{C}}

\newcommand {\shHom} {\mathcal{H}\!\text{\textit{om}}}
\newcommand {\shI}  {\mathcal{I}}

\newcommand {\shK}  {\mathcal{K}}
\newcommand {\shL}  {\mathcal{L}}
\newcommand {\shM}  {\mathcal{M}}

\newcommand {\shN}  {\mathcal{N}}

\newcommand {\Aut}  {\operatorname{Aut}}

\newcommand {\C} {\mathscr{C}}

\newcommand {\codim} {\operatorname{codim}}

\newcommand {\coker} {\operatorname{coker}}
\newcommand {\cyc}  {\mathrm{cyc}}

\newcommand {\eps}  {\varepsilon}

\newcommand {\eqand} {\quad\text{and}\quad}

\newcommand {\hol}  {\mathrm{hol}}
\newcommand {\Hom}  {\operatorname{Hom}}

\newcommand {\id}  {\operatorname{id}}
\newcommand {\im}  {\operatorname{im}}

\renewcommand {\ker } {\operatorname{kern}}

\newcommand {\loc} {\mathrm{loc}}

\newcommand {\lra}  {\longrightarrow}

\renewcommand{\O}  {\mathcal{O}}

\newcommand {\quadand} {\quad\text{and}\quad}

\newcommand {\reg}  {\mathrm{reg}}

\newcommand {\sing} {\mathrm{sing}}

\def\mydate{\ifcase\month \or January\or February\or March\or
April\or May\or June\or July\or August\or September\or October\or 
November\or December\fi \space\number\day,\space\number\year}
\begin{document}

\title[Lefschetz fibrations]{On the holomorphicity of genus two
Lefschetz fibrations}

\author[Bernd Siebert]{Bernd Siebert$^*$}
\address{Mathematische Fakult\"at, Ruhr-Universit\"at,
    44780 Bochum}
\curraddr{Universit\'e Paris 7,
    Institut de Math\'ematiques
    --- Topologie et G\'eo\-m\'e\-trie Alg\'ebriques ---
    2 place Jussieu,
    75251 Paris Cedex 05}
\email{bernd.siebert@ruhr-uni-bochum.de}

\author[Gang Tian]{Gang Tian}
\address{Department of Mathematics, MIT,
Cambridge, MA 02139--4307}
\email{tian@math.mit.edu}
\thanks{$^*$supported by the Heisenberg program of the DFG}


\dedicatory{January 17, 2002}

\begin{abstract}
We prove that any genus-2 Lefschetz fibration without reducible
fibers and with ``transitive monodromy'' is holomorphic. The latter
condition comprises all cases where the number of singular fibers
$\mu\in 10\NN$ is not congruent to $0$ modulo $40$. An auxiliary
statement of independent interest is the holomorphicity of symplectic
surfaces in $S^2$-bundles over $S^2$, of relative degree $\le 7$ over
the base, and of symplectic surfaces in $\CC\PP^2$ of degree $\le
17$.
\end{abstract}

\maketitle
\vspace{2ex}

\tableofcontents

%===========================================================
\section*{Introduction}\noindent
A differentiable Lefschetz fibration is a differentiable surjection
$p:M\to S^2$ of a closed oriented manifold $M$ with only finitely
many critical points of the form $t\circ p (z_1,\ldots,z_n)=\sum
z_i^2$. Here $z_1,\ldots,z_n$ and $t$ are complex coordinates on $M$
and $S^2$ respectively that are compatible with the orientations. It
is natural to ask how far differentiable Lefschetz fibrations are
from holomorphic ones. This question becomes even more interesting in
view of Donaldson's result on the existence of symplectic Lefschetz
pencils on arbitrary symplectic manifolds \cite{donaldson}.
Conversely, by an observation of Gompf total spaces of differentiable
Lefschetz fibrations have a symplectic structure that is unique up to
isotopy. The study of differentiable Lefschetz fibrations is
therefore essentially equivalent to the study of symplectic
manifolds.

In dimension $4$ apparent invariants of a Lefschetz fibration are the
genus of the non-singular fibers and the number and types of
irreducible fibers. By the work of Gromov and McDuff \cite{mcduff}
any genus-$0$ Lefschetz fibration is in fact holomorphic. Likewise,
for genus-$1$ the topological classification of elliptic fibrations
by Moishezon and Livn\'e \cite{moishezon} implies holomorphicity in
all cases. Our main theorem generalizes these results to genus $2$.
For the statement of our result recall that the mapping class group
$\text{MC}_2$ is a $\ZZ/2\ZZ$-extension of the braid group
$B(S^2,6)$. We say that a genus-$2$ Lefschetz fibration has
\emph{transitive monodromy} if the monodromy acts transitively on the
set of strands.

\begin{Thm.A}\label{genus-2 Lefschetz fibrations}
Let $p:M\to S^2$ be a genus-$2$ differentiable Lefschetz fibration
with transitive monodromy. If all singular fibers are irreducible
then $p$ is isomorphic to a holomorphic Lefschetz fibration.
\end{Thm.A}

Note that the conclusion of the theorem becomes false if we allow
reducible fibers, see e.g.\ \cite{ozbagcistipsicz}. The authors
expect that a genus-$2$ Lefschetz fibration with $\mu$ singular
fibers, $t$ of which reducible, is holomorphic if $t\le c\cdot\mu$
for some universal constant $c$. In particular, any genus-$2$
Lefschetz fibration should become holomorphic after fiber sum with
sufficiently many copies of the rational genus-$2$ Lefschetz
fibration with $20$ irreducible singular fibers. This problem should
also be solvable by the method presented in this paper.

In \cite{st99} we reduced Theorem~A to a holomorphicity statement for
certain symplectic surfaces in $S^2$-bundles over $S^2$. The main
theorem therefore follows essentially from the following isotopy
result for symplectic surfaces in rational ruled symplectic
$4$-manifolds.

\begin{Thm.B}\label{isotopy in rational ruled}
Let $p:M\to S^2$ be an $S^2$-bundle and $\Sigma\subset M$ a connected
surface symplectic with respect to a symplectic form that is
isotopic to a K\"ahler form. If $\deg(p|_\sigma)\le 7$ then $\Sigma$
is symplectically isotopic to a holomorphic curve in $M$, for some
choice of complex structure on $M$.
\end{Thm.B}

\begin{remark}\label{wlog Sigma>0}
By Gromov-Witten theory there exist surfaces $H,F\subset
M$, homologous to a section with self-intersection $0$ or $1$ and a
fiber, respectively, with $\Sigma\cdot H\ge 0$, $\Sigma\cdot F\ge
0$. It follows that $c_1(M) \cdot\Sigma>0$ unless $\Sigma$ is
homologous to a negative section. In the latter case
Proposition~\ref{fibered J} produces an isotopy to a section with
negative self-intersection number. The result follows then by the
classification of $S^2$-bundles with section. We may therefore add
the positivity assumption $c_1(M) \cdot\Sigma>0$ to the hypothesis of
the theorem. The complex structure on $M$ may then be taken generic,
thus leading to $\CC\PP^2$ or the first Hirzebruch surface
$\mathbb{F}_1= \PP(\O_{\CC\PP^1}\oplus \O_{\CC\PP^1}(1))$.
\qed
\end{remark}

Disconnectedness of the branch locus causes a difficulty to extend our
theorem on genus-$2$ Lefschetz fibrations to all cases. One
possibility may be to employ braid-theoretic arguments to reduce to
the connected case. We will discuss this in a future paper.

A similar result holds for surfaces of low degree in $\CC\PP^2$.

\begin{Thm.C}\label{isotopy in P^2}
Any symplectic surface in $\CC\PP^2$ of degree $d\le 17$ is
symplectically isotopic to an algebraic curve.
\end{Thm.C}

For $d=1,2$ this theorem is due to Gromov \cite{gromov}, for $d=3$ to
Sikorav \cite{sikorav} and for $d\le 6$ to Shevchishin
\cite{shevchishin}.

Together with the classification of symplectic structures on
$S^2$-bundles over $S^2$ by McDuff, Lalonde, A.K.~Liu and T.J.~Li (see
\cite{lalondemcduff} and references therein) our results indeed imply
a stronger classification of sympletic submanifolds up to Hamiltonian
symplectomorphism. Here we wish to add only the simple observation
that a symplectic isotopy of symplectic submanifolds comes from a
family of Hamiltonian symplectomorphisms.

\begin{proposition}
Let $(M,\omega)$ be a symplectic $4$-manifold and assume that
$\Sigma_t\subset M$, $t\in [0,1]$ is a family of symplectic
submanifolds. Then there exists a family $\Psi_t$ of Hamiltonian
symplectomorphisms of $M$ with $\Psi_0=\id$ and
$\Sigma_t=\Psi_t(\Sigma_0)$ for every $t$.
\end{proposition}

\proof
At a $P\in\Sigma_{t_0}$ choose complex Darboux coordinates $z=x+iy$,
$w=u+iv$ with $w=0$ describing $\Sigma_{t_0}$. In particular,
$\omega=dx\wedge dy+du\wedge dv$. For $t$
close to $t_0$ let $f_t$, $g_t$ be the functions describing
$\Sigma_t$ as graph $w=f_t(z)+i g_t(z)$. Define
$$
H_t= -(\partial_t g_t)\cdot (u-f_t) +(\partial_t f_t)\cdot(v-g_t).
$$
Then for every fixed $t$
$$
dH_t= -(u-f_t)\partial_t (dg_t) +(v-g_t)\partial_t (df_t)
-(\partial_t g_t)du+(\partial_t f_t)dv.
$$
Thus along $\Sigma_t$ 
$$
dH_t= -(\partial_t g_t)du+(\partial_t f_t)dv =\omega\neg\big((\partial_t
f_t)\partial_u+(\partial_t g_t) \partial_v\big).
$$
The Hamiltonian vector field belonging to $H_t$ thus induces the given
deformation of $\Sigma_t$.

To globalize patch the functions $H_t$ constructed locally
over $\Sigma_{t_0}$ by a partition of unity on $\Sigma_{t_0}$. As
$H_t$ vanishes along $\Sigma_t$, at time $t$ the associated Hamiltonian
vector field along $\Sigma_t$ remains unchanged. Extend $H_t$
to all of $M$ arbitrarily. Finally extend the construction to all
$t\in[0,1]$ by a partition of unity argument in~$t$.
\qed
\medskip

\emph{Conventions:}\ \ 
We endow complex manifolds such as $\CC\PP^n$ or $\mathbb{F}_1$ with
their integrable complex structures, when viewed as almost complex
manifolds. A map $F:(M,J_M)\to (N,J_N)$ of almost complex manifolds
is \emph{pseudo-holomorphic} if $DF\circ J_M= J_N\circ DF$. A
\emph{pseudo-holomorphic curve} $C$ in $(M,J)$ is the image of a
pseudo-holomorphic map $\varphi: (\Sigma,j)\to (M,J)$ with $\Sigma$ a
not necessarily connected Riemann surface. If $\Sigma$ may be chosen
connected then $C$ is \emph{irreducible} and its \emph{genus} $g(C)$
is the genus of $\Sigma$ for the generically injective $\varphi$. If
$g(C)=0$ then $C$ is \emph{rational}.

A \emph{$J$-holomorphic $2$-cycle} in an almost complex manifold
$(M,J)$ is a locally finite formal linear combination $C=\sum_a m_a
C_a$ where $m_a\in\ZZ$ and $C_a\subset M$ is a $J$-holomorphic curve.
The support $\bigcup_a C_a$ of $C$ will be denoted $|C|$. The subset
of singular and regular points of $|C|$ are denoted $|C|_\sing$ and
$|C|_\reg$ respectively. If all $m_a=1$ the cycle is \emph{reduced}.
We identify such $C$ with their associated pseudo-holomorphic curve
$|C|$. A \emph{smoothing} of a pseudo-holomorphic cycle $C$ is a
sequence $\{C_n\}$ of smooth pseudo-holomorphic cycles with $C_n\to
C$ in the $\C^0$-topology, see Section~\ref{cycle topology}. By abuse
of notation we often just speak of a smoothing $C^\dagger$ of $C$
meaning $C^\dagger=C_n$ with $n\gg0$ as needed. 

For an almost complex manifold $\Lambda^{0,1}$ denotes the bundle of
$(0,1)$-forms. \emph{Complex coordinates} on an even-dimensional,
oriented manifold $M$ are the components of an oriented chart
$M\supset U\to \CC^n$. Throughout the paper we fix some $0<\alpha<1$.
Almost complex structures will be of class $\C^l$ for some
sufficiently large integer $l$ unless otherwise mentioned. The unit
disk in $\CC$ is denoted $\Delta$. If $S$ is a finite set then $\#S$
is its cardinality. We measure distances on $M$ with respect to any
Riemannian metric, chosen once and for all. The symbol $\sim$ denotes
homological equivalence. An \emph{exceptional sphere} in an oriented
manifold is an embedded $2$-sphere with self-intersection number
$-1$.

%===========================================================
\section{Pseudo-holomorphic $S^2$-bundles}\label{sect1}\noindent
In our proof of the isotopy theorems it will be crucial to reduce to
a fibered situation. In Sections~\ref{sect1}, \ref{sect2} and
\ref{sect4} we introduce the notation and some of the tools that we
have at disposal in this case.

\begin{definition}
Let $p:M\to B$ be a smooth $S^2$-fiber bundle. If $M=(M,\omega)$ is a
symplectic manifold and all fibers $p^{-1}(b)$ are symplectic we
speak of a \emph{symplectic $S^2$-bundle}. If $M=(M,J)$ and $B=(B,j)$
are almost complex manifolds and $p$ is pseudo-holomorphic we speak
of a \emph{pseudo-holomorphic} $S^2$-bundle. If both preceding
instances apply and $\omega$ tames $J$ then $p:(M,\omega,J)\to (B,j)$
is a \emph{symplectic pseudo-holomorphic $S^2$-bundle}.
\end{definition}

In the sequel we will only consider the case $B=\CC\PP^1$. Then $M\to
\CC\PP^1$ is differentiably isomorphic to one of the holomorphic
$\CC\PP^1$-bundles $\CC\PP^1\times\CC\PP^1\to\CC\PP^1$ or
$\mathbb{F}_1\to\CC\PP^1$.

Any almost complex structure making a symplectic fiber bundle over a
symplectic base pseudo-holomorphic is tamed by some symplectic form.
To simplify computations we restrict ourselves to dimension~$4$.

\begin{proposition}\label{always tamed}
Let $(M,\omega)$ be a closed symplectic $4$-manifold and $p: M\to B$ a
smooth fiber bundle with all fibers symplectic. Then for any
symplectic form $\omega_B$ on $B$ and any almost complex structure $J$
on $M$ making the fibers of $p$ pseudo-holomorphic, $\omega_k:=
\omega+k\, p^*(\omega_B)$ tames $J$ for $k\gg0$.
\end{proposition}

\proof
Since tamedness is an open condition and $M$ is compact it suffices to
verify the claim at one point $P\in M$. Write $F=p^{-1}(p(P))$. Choose
a frame $\partial_u,\partial_v$ for $T_PF$ with
$$
J(\partial_u)=\partial_v,\quad
\omega(\partial_u,\partial_v)=1.
$$
Similarly let $\partial_x,\partial_y$ be a frame for the
$\omega$-perpendicular plane $(T_PF)^\perp\subset T_PM$ with
$$
J(\partial_x)=\partial_y+\lambda\partial_u+\mu\partial_v,\quad
\omega (\partial_x,\partial_y)=1
$$
for some $\lambda,\mu\in\RR$. By rescaling $\omega_B$ we may
also assume $(p^*\omega_B)(\partial_x,\partial_y) =1$. Replacing
$\partial_x,\partial_y$ by $\cos(t)\partial_x +\sin(t)\partial_y$,
$-\sin(t)\partial_x+\cos(t)\partial_y$, $t\in[0, 2\pi]$, the
coefficients $\lambda=\lambda(t)$, $\mu=\mu(t)$ vary in a compact
set. It therefore suffices to check that for $k\gg 0$
$$
\frac{\omega_{k-1} \big(\partial_x+\alpha\partial_u+\beta \partial_v,
J(\partial_x+\alpha\partial_u+\beta \partial_v) \big)}{
k+\alpha^2+\beta^2} =
1+\frac{\alpha\mu-\beta\lambda}{k+\alpha^2+\beta^2}
$$
is positive for all $\alpha,\beta\in\RR$. This term is minimal for
$$
\textstyle
\alpha=-\sqrt{\frac{k}{1+(\lambda/\mu)^2} },\quad
\beta= \sqrt{\frac{k}{1+(\mu/\lambda)^2} },
$$
where the value is $1-\sqrt{ \frac{\lambda^2+\mu^2}{4k} }$. This
is positive for $k>(\lambda^2+\mu^2)/4$.
\qed
\medskip

Denote by $T^{0,1}_{M,J}\subset T^\CC_M$ the anti-holomorphic tangent
bundle of an almost complex manifold $(M,J)$. Consider a submersion
$p:(M,J)\to B$ of an almost complex $4$-manifold with all fibers
pseudo-holomorphic curves. Let $z=p^*(u),w$ be complex coordinates on
$M$ with $w$ fiberwise holomorphic. Then
$$
T^{0,1}_{M,J}=\langle \partial_{\bar z} +a\partial_z +b \partial_w,
\partial_{\bar w} \rangle
$$
for some complex-valued functions $a,b$. Clearly, $a$ vanishes
precisely when $p$ is pseudo-holomorphic for some almost complex
structure on $B$. The Nijenhuis tensor $N_J: T_M\otimes T_M\to T_M$,
defined by
$$
4N_J(X,Y)= [JX,JY]-[X,Y]-J[X,JY]-J[JX,Y],
$$
is antisymmetric and $J$-antilinear in each entry. In dimension $4$
it is therefore completely determined by its value on a pair of
vectors that do not belong to a proper $J$-invariant subspace. For
the complexified tensor it suffices to compute
\begin{eqnarray*}
\lefteqn{N^\CC_J(\partial_{\bar z} +a\partial_z +b \partial_w,
\partial_{\bar w})}\hspace{2cm}\\
&=&-\frac{1}{2}[\partial_{\bar z} +a\partial_z +b \partial_w,
\partial_{\bar w}]
+\frac{i}{2}J[\partial_{\bar z} +a\partial_z +b \partial_w,
\partial_{\bar w}]\\
&=& \frac{1}{2}(\partial_{\bar w} a)\big(\partial_z -i J\partial_z\big)
+(\partial_{\bar w} b)\partial_w.
\end{eqnarray*}
Since $\partial_z-iJ\partial_z$ and $\partial_w$ are linearly
independent we conclude:

\begin{lemma}\label{integrable cx str}
An almost complex structure $J$ on an open set $M\subset \CC^2$ with
$T^{0,1}_{M,J}=\langle \partial_{\bar z} +a\partial_z +b \partial_w,
\partial_{\bar w} \rangle$ is integrable if and only if
$\partial_{\bar w} a=\partial_{\bar w} b=0$.
\qed
\end{lemma}

\begin{example}
Let $T^{0,1}_{M,J}=\langle \partial_{\bar z} +w \partial_w,
\partial_{\bar w} \rangle$. Then $z$ and $we^{-\bar z}$ are
holomorphic coordinates on $M$.
\end{example}

The lemma gives a convenient characterization of integrable complex
structures in terms of the functions $a,b$ defining $T^{0,1}_{M,J}$.
To globalize we need a connection for $p$. The interesting case will
be $p$ pseudo-holomorphic or $a=0$, to which we restrict from now on.

\begin{lemma}\label{connection and almost cx strs.}
Let $p:M\to B$ be a submersion endowed with a connection $\nabla$ and
let $j$ be an almost complex structure on $B$. Then the set of almost
complex structures $J$ making
$$
p:(M,J)\lra (B,j)
$$
pseudo-holomorphic is in one-to-one correspondence with pairs
$(J_{M/B}, \beta)$ where
\begin{enumerate}
\item $J_{M/B}$ is an endomorphism of $T_{M/B}$ with $J_{M/B}^2=-\id$.
\item $\beta$ is a homomorphism $p^*(T_B)\to T_{M/B}$ that is complex
anti-linear with respect to $j$ and $J_{M/B}$:
$$
\beta(j(Z))=-J_{M/B}(\beta(Z)).
$$
\end{enumerate}
Identifying $T_M=T_{M/B}\oplus p^*(T_B)$ via $\nabla$ the
correspondence is
$$
J=\left( \begin{array}{cc} J_{M/B}& \beta\\
0& j\end{array}\right).
$$
\end{lemma}

\proof
The only point that might not be immediately clear is the equivalence
of $J^2=-\id$ with complex anti-linearity of $\beta$. This follows by
computing
$$
J^2= \left( \begin{array}{cc} J_{M/B}^2
& J_{M/B}\circ\beta +\beta\circ j\\
0&j^2 \end{array}\right).
$$
\vspace{-4ex}

\qed

\begin{lemma}\label{adapted coordinates}
Let $p:(M,J)\to (B,j)$ be a pseudo-holomorphic submersion, $\dim
M=4$, $\dim B=2$. Then locally in $M$ there exists a differentiable map
$$
\pi: M\lra \CC
$$
inducing a pseudo-holomorphic embedding $p^{-1}(Q)\to \CC$ for every
$Q\in B$.

Moreover, to any such $\pi$ let
$$
\beta: p^*(T^{0,1}_{B,j})\lra T^{1,0}_{M/B, J_{M/B}}
$$
be the homomorphism provided by Lemma~\ref{connection and almost cx
strs.} applied to the connection belonging to $\pi$. Let $w$ be the
pull-back by $\pi$ of the linear coordinate on $\CC$ and $u$ a
holomorphic coordinate on $B$. Then $z:=p^*(u)$ and $w$ are complex
coordinates on $M$, and
$$
\beta(\partial_{\bar u})=-2bi\partial_w,
$$
for the $\CC$-valued function $b$ on $M$ with
\begin{eqnarray}\label{T^{0,1}}
T^{0,1}_{M,J}=\langle \partial_{\bar w},
\partial_{\bar z} +b\partial_w\rangle.
\end{eqnarray}
\end{lemma}

\proof
Since $p$ is pseudo-holomorphic, $J$ induces a complex structure on
the fibers $p^{-1}(Q)$, varying smoothly with $Q\in B$. Hence locally
in $M$ there exists a $\CC$-valued function $w$ that fiberwise
restricts to a holomorphic coordinate. This defines the
trivialization $\pi$.

In the coordinates $z,w$ define $b$ via $\beta(\partial_{\bar u})= -2
b i\partial_w$. Then
$$
J(\partial_{\bar z})=-i\partial_{\bar z}-2bi\partial_w,
$$
so the projection of $\partial_{\bar z}$ onto $T^{0,1}_{M,J}$ is
$$
(\partial_{\bar z}+iJ(\partial_{\bar z}))/2=\partial_{\bar z}+
b\partial_w.
$$
\vspace{-5ex}

\qed
\medskip

The two lemmas also say how to define an almost complex structure
making a given $p:M\to B$ pseudo-holomorphic, when starting from a
complex structure on the base, a fiberwise conformal structure, and a
connection for~$p$.
\medskip

For the symplectic isotopy problem we can reduce to a fibered
situation by the following device.

\begin{proposition}\label{fibered J}
Let $p:(M,\omega)\to S^2$ be a symplectic $S^2$-bundle. Let
$\Sigma\subset M$ be a symplectic submanifold. Then there exists an
$\omega$-tamed almost complex structure $J$ on $M$ and a map
$p': (M,J)\to \CC\PP^1$ with the following properties.
\begin{enumerate}
\item $p'$ is isotopic to $p$.
\item $p'$ is pseudo-holomorphic.
\item $\Sigma$ is $J$-holomorphic.
\end{enumerate}
Moreover, if $\{\Sigma_t\}_t$ is a family of symplectic submanifolds
there exist families $\{p'_t\}_t$ and $\{J_t\}$ with the analogous
properties for every $t$.
\end{proposition}

\proof
We explained in \cite{st99}, Proposition~4.1 how to obtain a symplectic
$S^2$-bundle $p':M\to \CC\PP^1$, isotopic to $p$, so that all
critical points of the projection $\Sigma\to \CC\PP^1$ are simple and
positive. This means that near any critical point there exist complex
coordinates $z,w$ on $M$ with $z=(p')^*(u)$ for some holomorphic
coordinate $u$ on $\CC\PP^1$ and so that $\Sigma$ is the zero locus
of $z-w^2$. We may take these coordinates in such a way that $w=0$
defines a symplectic submanifold. This property will enter below when
we discuss tamedness.

Since the fibers of $p'$ are symplectic the $\omega$-perpendicular
complement to $T_{M/\CC\PP^1}$ in $T_M$ defines a subbundle mapping
isomorphically to $(p')^*(T_{\CC\PP^1})$. This defines a connection $\nabla$ for
$p'$. By changing $\nabla$ slightly near the critical points we may
assume that it agrees with the connection defined by the projections
$(z,w)\to w$.

The coordinate $w$ defines an almost complex structure along the
fibers of $p'$ near any critical point. Since at $(z,w)=(0,0)$ the
tangent space of $\Sigma$ agrees with $T_{M/\CC\PP^1}$, this almost
complex structure is tamed at the critical points with respect to the
restriction $\omega_{M/\CC\PP^1}$ of $\omega$ to the fibers. Choose a
complex structure $J_{M/\CC\PP^1}$ on $T_{M/\CC\PP^1}$ that is
$\omega_{M/\CC\PP^1}$-tamed and that restricts to this fiberwise
almost complex structure near the critical points. 

By Lemma~\ref{connection and almost cx strs.} it remains to
define an appropriate endomorphism
$$
\beta: (p')^*(T_{\CC\PP^1})\lra T_{M/\CC\PP^1}.
$$
By construction of $\nabla$ and the local form of $\Sigma$ we may put
$\beta\equiv0$ near the critical points. Away from the critical
points let $z=(p')^*(u)$ and $w$ be complex coordinates as in
Lemma~\ref{adapted coordinates}. Then $\Sigma$ is locally a graph
$w=f(z)$. This graph will be $J=J(\beta)$-holomorphic if and only if
$$
\partial_{\bar z}f=b(z,f(z)).
$$
Here $b$ is related to $\beta$ via $\beta(\partial_{\bar
u})=-2bi\partial_w$. Hence this defines $\beta$ along $\Sigma$ away
from the critical points. We want to extend $\beta$ to all of $M$
keeping an eye on tamedness. For non-zero $X+Y\in
(p')^*(T_{\CC\PP^1})\oplus T_{M/\CC\PP^1}$ the latter requires
$$
0<\omega(X+Y,J(X+Y))= \omega(X,j(X))+\omega(Y,J_{M/\CC\PP^1}(Y))
+\omega(Y,\beta(X)).
$$
Near the critical points we know that $\omega(X,j(X))>0$ because
$w=0$ defines a symplectic submanifold. Away from the critical points
$X$ and $j(X)$ lie in the $\omega$-perpendicular complement of a
symplectic submanifold and therefore $\omega(X,j(X))>0$ too. Possibly
after shrinking the neighbourhoods of the critical points above we
may therefore assume that tamedness holds for $\beta=0$. By
construction it also holds with the already defined $\beta$ along
$\Sigma$. Extend this $\beta$ differentiably to all of $M$
arbitrarily. Let $\chi_\eps:M\to [0,1]$ be a function with
$\chi_\eps|_\Sigma\equiv 1$ and with support contained in an
$\eps$-neighbourhood of $\Sigma$. Then for $\eps$ sufficiently small,
$\chi_\eps\cdot\beta$ does the job.

If $\Sigma$ varies in a family argue analogously with an additional
parameter~$t$.
\qed
\medskip

In the next section we will see some implications of the fibered
situation for the space of pseudo-holomorphic cycles.

%===========================================================
\section{Pseudo-holomorphic cycles on pseudo-holomorphic
$S^2$-bundles}\label{sect2}\noindent
One advantage of having $M$ fibered by pseudo-holomorphic curves is
that it allows to describe $J$-holomorphic cycles by Weierstra\ss\
polynomials, cf.\ \cite{st00}. Globally we are dealing with sections
of a relative symmetric product. This is the topic of the present
section.

Throughout $p:(M,J)\to B$ is a pseudo-holomorphic $S^2$-bundle. To study
$J$-holomorphic curves $C\subset M$ of degree $d$ over $B$ we
consider the $d$-fold relative symmetric product $M^{[d]}\to
B$ of $M$ over $B$. This is the quotient of the $d$-fold
fibered product $M^d_B:=M\times_B \ldots \times_B M$ by the permutation
action of the symmetric group $S_d$. Set-theoretically $M^{[d]}$
consists of $0$-cycles in the fibers of $p$ of length $d$.

\begin{proposition}\label{M^{[d]}}
There is a well-defined differentiable structure on $M^{[d]}$,
depending only on the fiberwise conformal structure on $M$ over $B$.
\end{proposition}

\proof
Let $\Phi: p^{-1}(U)\to \CC\PP^1$ be a local trivialization of $p$
that is compatible with the fiberwise conformal structure, see the
proof of Lemma~\ref{integrable complex structure} for existence. Let
$u$ be a complex coordinate on $U$. To define a chart near a
$0$-cycle $\sum m_a P_a$ choose $P\in\CC\PP^1\setminus\{\Phi(P_a)\}$
and a biholomorphism $w:\CC\PP^1 \setminus\{P\}\simeq \CC$. The
$d$-tuples with entries disjoint from $\Phi^{-1}(P)$ give an open
$S_d$-invariant subset
$$
U\times \CC^d\subset M^d_B.
$$
Now the ring of symmetric polynomials on $\CC^d$ is free. A set of
generators $\sigma_1,\ldots, \sigma_d$ together with $z=p^*(u)$
provide complex, fiberwise holomorphic coordinates on $(U\times
\CC^d)/S_d \simeq U\times\CC^d\subset M^{[d]}$.

Different choices lead to fiberwise biholomorphic transformations.
The differentiable structure is therefore well-defined.
\qed
\medskip

We emphasize that different choices of the fiberwise conformal
structure on $M$ over $B$ lead to different differentiable structures
on $M^{[d]}$. Note also that $M^d_B\lra M^{[d]}$ is a branched Galois
covering with covering group $S_d$. The branch locus is stratified
according to partitions of $d$, parametrizing cycles with the
corresponding multiplicities. The \emph{discriminant locus} is the
union of all lower-dimensional strata. The stratum belonging to a
partition $d=m_1+\ldots+m_1+\ldots+m_s +\ldots+ m_s$ with $m_1< m_2
< \ldots< m_s$ and $m_i$ occurring $d_i$-times is canonically
isomorphic to the complement of the discriminant locus in
$M^{[d_1]}\times_B\ldots\times_B M^{[d_s]}$. 

\begin{proposition}\label{cx structure on M^{[d]}}
There exists a unique continuous almost complex structure on
$M^{[d]}$ making the covering
$$
M^d_B\lra M^{[d]}
$$
pseudo-holomorphic.
\end{proposition}

\proof
It suffices to check the claim locally in $M^{[d]}$. Let
$w:U\times\CC\to \CC$ be a fiberwise holomorphic coordinate as in the
previous lemma. Let $z=p^*(u)$ and $b$ be as in Lemma~\ref{adapted
coordinates}, so
$$
T^{0,1}_{M,J}=\big\langle \partial_{\bar w},
\partial_{\bar z} +b(z,w)\partial_w \big\rangle.
$$
Let $w_i$ be the pull-back of $w$ by the $i$-th projection $M^d_B \to
M$. By the definition of the differentiable structure on $M^{[d]}$,
the $r$-th elementary symmetric polynomial $\sigma_r(w_1,\ldots,w_d)$
descends to a locally defined smooth function $\sigma_r$ on
$M^{[d]}$. The pull-back of $u$ to $M^{[d]}$, also denoted by $z$,
and the $\sigma_r$ provide local complex coordinates on $M^{[d]}$.
The almost complex structure on $M^d_B$ is
$$
T^{0,1}_{M^d_B}= \big\langle \partial_{\bar w_1},\ldots,
\partial_{\bar w_d}, \partial_{\bar z}+ b(z,w_1)\partial_{w_1}+\ldots+
b(z,w_d)\partial_{w_d}\big\rangle.
$$
The horizontal anti-holomorphic vector field
$$
\partial_{\bar z}+ b(z,w_1) \partial_{w_1}+\ldots+b(z,w_d)\partial_{w_d}
$$
is $S_d$-invariant, hence descends to a continuous vector field $Z$
on $M^{[d]}$. Together with the requirement that fiberwise the map
$M^d_B\to M^{[d]}$ be holomorphic, $Z$ determines the almost complex
structure on $M^{[d]}$. 
\qed

\begin{remark}
The horizontal vector field $Z$ in the lemma, hence the almost
complex structure on $M^{[d]}$, will generally only be of H\"older
class $\C^{0,\alpha'}$ in the fiber directions, for some
$\alpha'>0$,  see \cite{st99}. However, at $0$-cycles $\sum m_\mu
P_\mu$ with $J$ integrable near $\{P_\mu\}$ it will be smooth and
integrable as well. In fact, by Lemma~\ref{integrable cx str}
integrability is equivalent to holomorphicity of $b$ along the
fibers. This observation will be crucial below.
\end{remark}

\begin{proposition}\label{cycles-section}
There is an injective map from the space of $J$-holo\-mor\-phic
cycles on $M$ of degree $d$ over $B$ and without fiber components to
the space of $J_{M^{[d]}}$-holomorphic sections of $q: M^{[d]}\to B$.
A cycle $C=\sum m_a C_a$ maps to the section
$$
u\longmapsto \sum_a m_a C_a\cap p^{-1}(u),
$$
the intersection taken with multiplicities. The image of the subset
of reduced cycles are the sections with image not entirely lying in
the discriminant locus.
\end{proposition}

\proof
We may reduce to the local problem of cycles in $\Delta\times\CC$.
In this case the statement follows from \cite{st00}, Theorem~I. 
\qed

\begin{remark}
By using the stratification of $M^{[d]}$ by fibered products
$M^{[d_1]} \times_B\ldots\times_B M^{[d_k]}$ with $\sum d_i\le d$ it
is also possible to treat cycles with multiple components. In fact,
one can show that a pseudo-holomorphic section of $M^{[d]}$ has image
in exactly one stratum except at finitely many points. Now the almost
complex structure on a stratum agrees with the almost complex induced
from the factors. The claim thus follows from the proposition applied
to each factor. Because this result is inessential to what follows
details are left to the reader.
\end{remark}
\medskip

To study deformations of a $J$-holomorphic cycle it therefore
suffices to look at deformations of the associated section of
$M^{[d]}$. Essentially this is what we will do; but as we also have
to treat fiber components we describe our cycles by certain
polynomials with coefficients taking values in holomorphic line
bundles over $B$. We restrict ourselves to the case $B=\CC\PP^1$.

The description depends on the choice of an integrable complex
structure on $M$ fiberwise agreeing with $J$. Thus we assume now that
$p:(M,J_0)\to\CC\PP^1$ is a holomorphic $\CC\PP^1$-bundle. There
exist disjoint sections $S,H\subset M$ with $e:=H\cdot H\ge 0$. Then
$H\sim S+eF$ where $F$ is a fiber, and $S\cdot S=-e$. Denote the
holomorphic line bundles corresponding to $H,S$ by $L_H$ and $L_S$.
Let $s_0,s_1$ be holomorphic sections of $L_S,L_H$ respectively with
zero loci $S$ and $H$. We also choose an isomorphism $L_H\simeq
L_S\otimes p^*(L)^e$, where $L$ is the holomorphic line bundle on
$\CC\PP^1$ of degree $1$.

Note that if $H_d\subset M^{[d]}$ denotes the divisor of cycles
$\sum_a m_a P_a$ with $P_a\in H$ for some $a$ then
$$
M^{[d]}\setminus H_d= S^d(L^{-e})=\bigoplus_{\nu=1}^d L^{-\nu e}.
$$
In fact, $M\setminus H= L^{-e}$.

\begin{proposition}\label{cycles-coefficients}
Let $J$ be an almost complex structure on $M$ making
$p:M \to\CC\PP^1$ pseudo-holomorphic and assume $J=J_0$ near
$H$ and along the fibers of $p$.
\smallskip

\noindent
1)\ \ Let $C=\sum_a m_a C_a$ be a $J$-holomorphic $2$-cycle
homologous to $dH+kF$, $d>0$, and assume $H\not\subset |C|$. Let
$a_0$ be a holomorphic section of $L^{k+de}$ with zero locus
$p_*(H\cap C)$, with multiplicities.

Then there are unique continuous sections $a_r$ of $L^{k+(d-r)e}$,
$r=1,\ldots, d$, so that $C$ is the zero locus of
$$
p^*(a_0) s_0^d+ p^*(a_1) s_0^{d-1} s_1+\ldots
+p^*(a_d) s_1^d,
$$
as cycle.
\smallskip

\noindent
2)\ \ There exist H\"older continuous maps
$$
\beta_r:\bigoplus_{\nu=1}^d L^{-\nu e}\lra L^{-re}\otimes
\Lambda^{0,1}_{\CC\PP^1},\quad r=1,\ldots,d,
$$
so that a local section $s_0^d+p^*(\alpha_1)s_0^{d-1} s_1+\ldots+
p^*(\alpha_d) s_1^d=0$ of $M^{[d]}\setminus H_d$ is
$J_{M^{[d]}}$-holomorphic if and only if
$$
\bar\partial \alpha_r = \beta_r(\alpha_1,\ldots,\alpha_d),\quad r=1,\ldots,d.
$$

\noindent
3)\ \ Let $C$ be a $J$-holomorphic $2$-cycle homologous to $dH+kF$
and with $H\not\subset |C|$. Decompose $C=\bar C+ \sum m_a F_a$ with
the second term containing all fiber components. Assume that $J=J_0$
also near
$$
p^{-1}\big(p(|\bar C|\cap H)\cap p(|\bar C|\cap S) \big)\cup
{\textstyle\bigcup_a F_a}.
$$
Let $a^0_r$ be sections of $L^{k+(d-r)e}$ associated to $C$ according
to (1). Then there exists a neighbourhood $\mathfrak{D}\subset
\bigoplus_{r=0}^d L^{k+(d-r) e}$ of the graph of
$(a^0_0,\ldots,a^0_d)$ and H\"older continuous maps
$$
b_r: \mathfrak{D}\lra L^{k+(d-r)e}\otimes \Lambda^{0,1}_{\CC\PP^1},
\quad r=1,\ldots,d,
$$
so that a section $(a_0,\ldots,a_d)$ of $\mathfrak{D}\to\CC\PP^1$
with $a_0$ holomorphic defines a $J$-holomorphic cycle if and only if
\begin{eqnarray}\label{PDE}
\bar\partial a_r = b_r(a_0,\ldots,a_d),\quad r=1,\ldots,d.
\end{eqnarray}
Conversely, any solution of {\rm (\ref{PDE})}\! with
$\delta(a_0,\ldots,a_d)\not\equiv 0$ corresponds to a $J$-holomorphic
cycle without multiple components. Here $\delta$ is the discriminant.

Moreover, if $J$ is integrable near $|C|$ then the $b_r$ are smooth
near the corresponding points of $\mathfrak{D}$.
\end{proposition}

\proof
1)\ \ Assume first that $a=1$ and $m_1=1$. Then either $C$ is a fiber
and $p^*(a_0)$ is the defining polynomial; or $C$ defines a section
of $M^{[d]}$ as in Proposition~\ref{cycles-section}. Any $0$-cycle of
length $d$ on $p^{-1}(Q)\simeq \CC\PP^1$ is the zero locus of a
section of $\O_{\CC\PP^1}(d)$ that is unique up to rescaling. The
restrictions of $s_0^r s_1^{d-r}$ to any fiber form a basis for the
space of global sections of $\O_{\CC\PP^1}(d)$. Hence, after choice
of $a_0$ the $a_r$ are determined uniquely for $r=1,\ldots,d$ away
from the zero locus of $a_0$.  If $a_0(Q)=0$ choose a neighbourhood
$U$ of $Q$ so that $C|_{p^{-1}(U)}= C'+C''$ with $|C'|\cap
S=\emptyset$, $|C''|\cap H=\emptyset$. By the same argument as before
we have unique Weierstra\ss\ polynomials of the form
\begin{eqnarray*}
&p^*(a'_0) s_0^{d'}+\ldots+ p^*(a'_{d'-1}) s_0s_1^{d'-1}+ s_1^{d'},\\
&s_0^{d''}+p^*(a''_1) s_0^{d''-1}s_1+\ldots+ p^*(a''_{d''})
s_1^{d''}
\end{eqnarray*}
defining $C'$ and $C''$ respectively. Multiplying produces a
polynomial defining $C$. The first coefficient $a'_0$ vanishes to the
same order at $Q$ as $a_0$. In fact, this order equals the
intersection index of $C'$ and $C$ with $H$ respectively. This shows
$a_0=a'_0\cdot e$ for some holomorphic function $e$ on $U$ with
$e(Q)\neq 0$. Therefore $a_1,\ldots,a_d$ extend continuously over
$Q$.

In the general case let $F_a$ be the polynomial just obtained for
$C=C_a$. Put
$$
F_{(a_0,\ldots, a_d)} = \prod_a F_a^{m_a}.
$$
The coefficient of $s_0^d$ has the same zero locus as $p^*(a_0)$; so after
rescaling by a constant, $F_{(a_0,\ldots, a_d)}$ has the desired form.
\smallskip

\noindent
2)\ \ Since $J$ and $J_0$ agree fiberwise and both make $p$
pseudo-holomorphic, the homomorphism $J-J_0$ factors over $p^*
T_{\CC\PP^1}$ and takes values in $T_{M/\CC\PP^1}$. Let $\beta$ be
the section of $T^{1,0}_{M/\CC\PP^1,J}\otimes
p^*\Lambda^{0,1}_{\CC\PP^1}$ thus defined. Locally $\beta$ is nothing
but the homomorphism obtained by applying Lemma~\ref{adapted
coordinates}. to a local $J_0$-holomorphic trivialization of $M$.
Because $M\setminus H= L^{-e}$ there is a canonical isomorphism
$$
T^{1,0}_{M/\CC\PP^1, J_0}\big|_{M\setminus H} \simeq p^*(L^{-e}).
$$
This isomorphism understood we obtain an $S_d$-invariant map
$$
(w_1,\ldots,w_d)\longmapsto (-1)^r\frac{i}{2} \sum_\nu \sigma_{r-1}
(w_1,\ldots, \widehat w_\nu, \ldots, w_d) \otimes \beta(z,w_\nu)
$$
from $(L^{-e})^{\oplus d}$ to $L^{-re}\otimes \Lambda^{0,1}_{\CC\PP^1}$.
Define $\beta_r(\alpha_1,\ldots,\alpha_d)$ as the induced map from
$S^d(L^{-e})=M^{[d]}\setminus H_d$. The claim on pseudo-holomorphic
sections of $M^{[d]}\setminus H_d$ is clear from the definition of
$J_{M^{[d]}}$ in Proposition~\ref{cx structure on M^{[d]}} and the
description of $\beta$ in Lemma~\ref{adapted coordinates}.

H\"older continuity of the $\beta_r$ follows from the
local consideration in \cite{st99}.
\smallskip

\noindent
3)\ \ Let $U$ be a neighbourhood of $p^{-1}\big(p(|\bar C|\cap H)\cap
p(|\bar C|\cap S) \big) \cup \bigcup F_a$ with $J=J_0$ on
$p^{-1}(U)$. Over $Q\in\CC\PP^1$ define $\mathfrak{D}$ as those tuples
$(a_0,\ldots,a_d)$ with
$$
a_0=0\quad\Longrightarrow\quad a_d\neq 0 \quad \text{or}\quad Q\in U.
$$
If $a_s=0$ for $s=0,\ldots,m-1$ and $a_m\neq0$ define
$$
b_r(a_0,\ldots,a_d):= a_m\cdot
\beta^{d-m}_{r-m}(a_{m+1}/a_m,\ldots,a_d/a_m),
$$
where $\beta^{d'}_r$ is $\beta_r$ from (2) for $d=d'$. We also
put $b_r(0,\ldots,0)=0$. We claim that the $b_r$ are continuous. Over
$U$ this is clear as all terms vanish. 

Let $w$ be a complex coordinate on $M$ defining a local
$J_0$-holomorphic trivialization of $M\setminus H\to \CC\PP^1$. Let
$w_1,\ldots,w_d$ be the induced coordinate functions on
$M^d_{\CC\PP^1}$ and $b$ the function encoding $J$. Pulling back
$b_r$ via $M^d_{\CC\PP^1}\to M^{[d]}$ gives
\begin{eqnarray}\label{b_r}
a_0\cdot \sum_{\nu=1}^d \sigma_{r-1}(w_1,\ldots,\widehat
w_\nu,\ldots,w_d) b(z,w_\nu).
\end{eqnarray}
It remains to show that if $\{\lambda^{(n)}_1, \ldots,
\lambda_d^{(n)}\}_n$ and $\{a_0^{(n)}\}_n$ are sequences with
$a_r^{(n)}:=a_0^{(n)} \sigma_r (\lambda_1^{(n)},\ldots,
\lambda_d^{(n)})$ converging to $(0,\ldots,0,a_m,\ldots,a_d)$ with
$a_m\neq0$, $a_d\neq0$, then (\ref{b_r}) converges towards $a_m\cdot
\beta^{d-m}_{r-m} (a_{m+1}/a_m, \ldots,a_d/a_m)$. After reordering we
may assume that $\lambda_1^{(n)},\ldots,\lambda_m^{(n)}$ correspond
to the $m$ points converging to $H$. By hypothesis $b(z,w)=0$ for
$|w|\gg0$. Moreover, since $a^{(n)}_d$ converges with non-zero limit
and all $a_r^{(n)}$ are bounded, the $\lambda^{(n)}_\nu$ stay
uniformly bounded away from $0$. Hence for any subset
$I\subset\{1,\ldots,d\}$
$$
a^{(n)}_0\prod_{\nu\in I}\lambda^{(n)}_\nu
$$
converges. The limit is $0$ if $\{1,\ldots,m\}\not\subset I$.
Evaluating expression (\ref{b_r}) at $w_\nu=\lambda^{(n)}_\nu$ and
taking the limit gives
\begin{eqnarray*}
&&\lim_{n\to\infty} \Big(a^{(n)}_0\cdot \sum_{\nu=m+1}^d \sigma_{r-1}
(\lambda^{(n)}_1,\ldots,\widehat
\lambda^{(n)}_\nu,\ldots,\lambda^{(n)}_d) \cdot
b(z,\lambda^{(n)}_\nu)\Big)\\
&&=\ a_m\cdot\lim_{n\to\infty} \sum_{\nu=m+1}^d \sigma_{r-m-1}
(\lambda^{(n)}_{m+1},\ldots,\widehat
\lambda^{(n)}_\nu,\ldots,\lambda^{(n)}_d)  \cdot
b(z,\lambda^{(n)}_\nu)\\
&&=\ a_m \cdot \beta^{d-m}_{r-m}(a_{m+1}/a_m,\ldots,a_d/a_m),
\end{eqnarray*}
as had to be shown.

The expression for $b_r$ also shows that the local equation for
pseudo-holo\-mor\-phicity of a section $\sigma_r=a_r(z)/a_0(z)$ of
$M^{[d]}\setminus H_d$ is
$$
\partial_{\bar z} a_r(z)= a_0 \beta_r(a_1,\ldots,a_d)
=b_r(a_0,\ldots,a_d).
$$
This extends over the zeros of $a_0$. The converse follows from the
local situation already discussed at length in \cite{st00}. 
\smallskip

Finally we discuss regularity of the $b_r$. The partial derivatives
of $b_r$ in the $z$-direction lead to expressions of the same form as
$b_r$ with $b(z,w_\nu)$ replaced by $\nabla^k_z b(z,w_\nu)$. These
are continuous by the argument in (2). If $J$ is integrable near
$|C|$ then $b$ is holomorphic there along the fibers of $p$. Hence
the $b_r$ and its derivatives in the $z$-direction are continuous and
fiberwise holomorphic. Uniform boundedness thus implies the desired
estimates on higher mixed derivatives. 
\qed

\begin{remark}
It is instructive to compare the linearizations of the equations
characterizing $J$-holomorphic cycles of the coordinate dependent
description in this proposition and the intrinsic one in
Proposition~\ref{cycles-section}. We have to assume that $C$ has no
fiber components. Let $\sigma$ be the section of $q:M^{[d]}\to \CC\PP^1$
associated to $C$ by Proposition~\ref{cycles-section}. There is a PDE
acting on sections of $\sigma^*(T_{M^{[d]}/\CC\PP^1})$ governing
(pseudo-) holomorphic deformations of $\sigma$. For the integrable
complex structure this is simply the $\bar\partial$-equation.
There is a well-known exact sequence
$$
0\lra \underline{\CC}\lra \bigoplus_{r=0}^d L^{k+(d-r)e}
\lra \sigma^*(T_{M^{[d]}/\CC\PP^1}\big) \lra 0,
$$
describing the pull-back of the relative tangent bundle. The $\bar
\partial_J$-equation giving $J$-holomorphic deformations of
$\sigma$ acts on the latter bundle. On the other hand, the middle
term exhibits variations of the coefficients $a_0,\ldots,a_d$. The
constant bundle on the left deals with rescalings.
\end{remark}
\bigskip

The final result of this section characterizes certain smooth cycles.

\begin{proposition}\label{smooth cycles}
In the situation of Proposition~\ref{cycles-section} let $\sigma$ be
a differentiable section of $M^{[d]}\to S^2$ intersecting the discriminant
divisor transversally. Then the 2-cycle $C$ belonging to $\sigma$ is
a submanifold and the projection $C\to S^2$ is a branched cover with
only simple branch points. Moreover, $C$ varies differentiably under
$\mathscr{C}^1$-small variations of $\sigma$.
\end{proposition}

\proof
Away from points of intersection with the discriminant divisor the
symmetrization map $M^d_B\to M^{[d]}$ is locally a diffeomorphism,
and the result is clear. Moreover, the discriminant divisor is smooth
only at points $\sum m_\mu P_\mu$ with $\sum m_\mu=d-1$; this is the
locus where exactly two points come together. We may hence assume
$m_1=2$ and $m_a=1$ for $a>1$. At $\sum m_\mu P_\mu$ the $d-2$
coordinates $w_\mu$ at $P_\mu$, $\mu>2$, and $w_1 +w_2$, $w_1w_2$
descend to complex coordinates on $S^d(\CC\PP^1)$. Similarly, the
variation of the $P_\mu$ for $\mu>2$ lead only to multiplication of
$\delta(a_0,\ldots,a_d)$ by a smooth function without zeros. It
therefore suffices to discuss the case $d=2$. Then $C$ is the zero
locus
$$
a_0(z)w^2+a_1(z)w+a_2(z)=0.
$$
The assumption says that, say, $z=0$ is a simple zero of
$$
\delta(\alpha_1,\alpha_2)=\alpha_1^2-4\alpha_2.
$$
By assumption there exists a function $h(z)$ with $h(0)\neq0$ and
$\delta(\alpha_1,\alpha_2) =h^2(z)\cdot z$. Replacing $w$ by
$u=2h^{-1}(w+\frac{\alpha_1}{2})$ brings $C$ into standard from
$u^2-z=0$. Hence $C$ is smooth and the projection to $z$ has a simple
branch point over $z=0$. The same argument is valid for small
deformations of~$\sigma$.
\qed

%===========================================================
\section{The $\C^0$-topology on the space of
pseudo-holomorphic cycles}\label{cycle topology}\noindent
This section contains a discussion of the topology on the space of
pseudo-holomorphic cycles, which we denote
$\operatorname{Cyc}_{\mathrm{pshol}}(M)$ throughout. Let $\shC(M)$ be
the space of pseudo-holomorphic stable maps. An element of $\shC(M)$
is an isomorphism class of pseudo-holomorphic maps $\varphi:\Sigma\to
M$ where $\Sigma$ is a nodal Riemann surface, with the property that
there are no infinitesimal biholomorphisms of $\Sigma$ compatible with
$\varphi$. The $\C^0$-topology on $\shC(M)$ is generated by open sets
$U_{V,\eps}$ defined for $\eps>0$ and $V$ a neighbourhood of
$\Sigma_\sing$ as follows. To compare $\psi:\Sigma'\to M$ with
$\varphi$ consider maps $\kappa: \Sigma'\to\Sigma$ that are a
diffeomorphism away from $\Sigma_\sing$ and that over a branch of
$\Sigma$ at a node have the form
$$
\big\{z\in\Delta\,\big|\, |z|>\tau\big\}\lra \Delta,\quad
r e^{i\phi}\longmapsto \frac{r-\tau}{1-\tau}e^{i\phi}
$$
for some $0\le\tau<1$. Then $\psi:\Sigma'\to M$ belongs to
$U_{V,\eps}$ if such a $\kappa$ exists with maximal dilation over
$\Sigma\setminus V$ less than $\eps$ and with
$$
d_M\big(\psi(z),\varphi(\kappa(z))\big)<\eps\quad
\text{for all $z$.}
$$
Recall that the dilation measures the failure of a map between
Riemann surfaces to be holomorphic. Note also that an intrinsic measure for
the size of the neighbourhood $V$ of the singular set on non-contracted
components is the diameter of $\varphi(V)$ in $M$; on contracted
components one may take the smallest $\eps$ with $V$ contained in the
\emph{$\eps$-thin part}. The latter consists of endpoints of
loops around the singular points of length $<\eps$ in the Poincar\'e
metric.

For a fixed almost complex structure of class
$\mathscr{C}^{l,\alpha}$, $\C^0$-convergence of
pseudo-holo\-mor\-phic stable maps implies
$\C^{l+1,\alpha}$-convergence away from the singular points of the
limit. If one wants convergence of derivatives away from the
singularities for varying $J$ one needs $\C^{0,\alpha}$-convergence
of $J$ for some $\alpha>0$. We will impose this condition separately
each time we need it.

The $\C^0$-topology on $\shC(M)$ induces a topology on
$\operatorname{Cyc}_{ \mathrm{pshol}}(M)$ via the map
$$
\shC(M)\lra \operatorname{Cyc}_{\mathrm{pshol}}(M),\quad
\big(\varphi:C={\textstyle\bigcup_a} C_a\to M\big)
\longmapsto \sum_a m_a\varphi(C_a).
$$
Here $m_a$ is the covering degree of $C_a\to \varphi(C_a)$. From this
point of view the compactness theorem for $\operatorname{Cyc}_{
\mathrm{pshol}}(M)$ follows immediately from the version for stable
maps. We call this topology on $\operatorname{Cyc}_{
\mathrm{pshol}}(M)$ the \emph{$\C^0$-topology}. 

Alternatively, one may view $\operatorname{Cyc}_{ \mathrm{pshol}}(M)$
as a closed subset of the space of currents on $M$, or of the space
of measures on $M$. We will not use this point of view here.
\bigskip

Next we turn to a semi-continuity property of pseudo-holomorphic
cycles in the $\C^0$-topology. For a pseudo-holomorphic curve
singularity $(C,P)$ in a $4$-dimensional almost complex manifold $M$
define $\delta(C,P)$ as the virtual number of double points. This is
the number of nodes of the image of a small, general, $J$-holomorphic
deformation of the parametrization map from a union of unit disks to
$M$ belonging to $(C,P)$. This number occurs in the genus formula. If
$C=\bigcup_{a=1}^d C_a$ is the decomposition of a pseudo-holomorphic
curve into irreducible components, the genus formula says
\begin{eqnarray}\label{genus formula}
\sum_{a=1}^d g(C_a)=\frac{C\cdot C-c_1(M)\cdot C}{2}+d
-\sum_{P\in C_\sing} \delta(C,P).
\end{eqnarray}
We emphasize that in this formula $C$ has no multiple components.
For a proof perform a small, general pseudo-holomorphic
deformation $\varphi_a:\Sigma_a\to M$ of the pseudo-holo\-mor\-phic
maps with image $C_a$. This is possible by changing $J$ slightly away
from $C_\sing$. The result is a $J'$-holomorphic nodal curve for some
small perturbation $J'$ of $J$. The degree of the complex line bundle
$\varphi_a^*(T_M)/T_{\Sigma_a}$ equals $C_a\cdot C_a$ minus the
number of double points of $C_a$. This expresses the genus of
$\Sigma_a$ in terms of $C_a\cdot C_a$, $c_1(M)\cdot C_a$ and
$\sum_{P\in (C_a)_\sing} \delta(C_a,P)$. Sum over $a$ and adjust by
the intersections of $C_a$ with $C_{a'}$ for $a\neq a'$ to deduce
(\ref{genus formula}). 

As a measure for how far a pseudo-holomorphic cycle differs from a
cycle with only ordinary double points we introduce
$$
\delta(C):=\sum_{P\in |C|_\sing}  \big(\delta(|C|,P) -1\big).
$$
Similarly, as a measure for non-reducedness of a pseudo-holomorphic
cycle $C=\sum_a m_a C_a$ put
$$
m(C):=\sum_a (m_a-1).
$$
So $\delta(C)=0$ if and only if $|C|$ has only ordinary double
points and $m(C)=0$ if and only if $C$ has no multiple components.

The definition of the $\C^0$-topology on the space of
pseudo-holomorphic cycles implies semi-continuity of the pair
$(m,\delta)$.

\begin{lemma}\label{semicontinuity of m,delta}
Let $(M,J)$ be a $4$-dimensional almost complex manifold with $J$
tamed by some symplectic form. Let $C_n\subset M$ be
$J_n$-holomorphic cycles with $J_n\to J$ in $\C^0$ and in
$\C^{0,\alpha}$ away from a set not containing any closed
pseudo-holomorphic curves, and assume $C_n\to C_\infty$ in the
$\C^0$-topology.

Then $m(C_n)\le m(C_\infty)$ for $n\gg0$, and if $m(C_n)=
m(C_\infty)$ for all $n$ then $\delta(C_n)\le \delta(C_\infty)$.
Moreover, if also $\delta(C_n) = \delta(C_\infty)$ for all $n$ then
for $n\gg 0$ there is a bijection between the irreducible components
of $|C_n|$ and of $|C_\infty|$ respecting the genera.
\end{lemma}

\proof
By the definition of the cycle topology, for $n\gg0$ each component
of $C_\infty$ deforms to parts of some component of $C_n$. This sets
up a surjective multi-valued map $\Delta$ from the set of irreducible
components of $C_\infty$ to the set of irreducible components of
$C_n$. The claim on semi-continuity of $m$ follows once we show that
the sum of the multiplicities of the components $C_{n,i}\in
\Delta(C_{\infty,a})$ does not exceed the multiplicity of
$C_{\infty,a}$.

By the compactness theorem we may assume that the $C_n$ lift to a
converging sequence of stable maps $\varphi_n:\Sigma_n\to M$. Let
$\varphi_\infty: \Sigma_\infty\to M$ be the limit. This is a stable
map lifting $C_\infty$. For a component $C_{\infty,a}$ of $C_\infty$
of multiplicity $m_a$ choose a point $P\in C_{\infty,a}$ in the part
of $\C^{0,\alpha}$-convergence of the $J_n$ and away from the
critical values of $\varphi_\infty$. Let $H\subset M$ be a local
submanifold of real codimension $2$ intersecting $|C|$ transversally
and positively in $P$. As $\C^{0,\alpha}$-convergence of almost
complex structures implies convergence of tangent spaces away from
the critical values, $H$ is transverse to $C_n$ for $n\gg0$ with all
intersections positive. Now any component of $C_n$ with a part
degenerating to $C_{\infty,a}$ intersects $H$, and $H\cdot C_n$
gives the multiplicity of $C_{\infty,a}$ in $C_\infty$. The claimed
semi-continuity of multiplicities thus follows from the deformation
invariance of intersection numbers.

The argument also shows that equality $m(C_\infty)=m(C_n)$ can only
hold if $\Delta$ induces a bijection between the non-reduced
irreducible components of $C_n$ and $C_\infty$ respecting the
multiplicities. This implies convergence $|C_n|\to |C_\infty|$, so we
may henceforth assume $C_n$ and $C_\infty$ to be reduced, and
$\varphi_n$ to be injective. Note that $\varphi_\infty$ may contract
some irreducible components of $\Sigma_\infty$. In any case, $\sum_a
g(C_{\infty,a})$ is the sum of the genera of the non-contracted
irreducible components of $\Sigma_\infty$, and it is not larger than
the respective sum for $\Sigma_n$. The latter equals $\sum_i
g(C_{n,i})$ if $C_n=\bigcup_i C_{n,i}$. By the genus formula
(\ref{genus formula}) we conclude
$$
\sum_{Q\in (C_n)_\sing}\delta(C_n,Q)\le
\sum_{P\in (C_\infty)_\sing}\delta(C_\infty,P),
$$
for $C_n\sim C_\infty$ in homology and because the number of
irreducible components is semicontinuous by surjectivity of $\Delta$.
Moreover, the number of singular points can at most decrease in the
limit. This shows $\delta(C_n)\le \delta(C_\infty)$. If equality
holds, there is a bijection between the singular points of $|C_n|$
and $|C_\infty|$ respecting the number of virtual double points. The
genus formula then shows that $\Delta$ respects the genera of the
irreducible components. 
\qed
\bigskip

In the fibered situation of Proposition~\ref{cycles-section}
convergence in $\operatorname{Cyc}_{\mathrm{pshol}}(M)$ in the
$\C^0$-topology implies convergence of the section of $M^{[d]}$:

\begin{proposition}\label{section convergence}
Let $p:M\to B$ be an $S^2$-bundle. For every $n$ let $C_n$ be a
pseudo-holomorphic curve of degree $d$ over $B$ for some almost
complex structure making $p$ pseudo-holomorphic. Assume that $C_n\to
C$ in the $\C^0$-topology and that $C$ contains no fiber components.
Let $\sigma_n$ and $\sigma$ be the sections of $M^{[d]}\to B$
corresponding to $C_n$ and $C$, respectively, according to
Proposition~\ref{cycles-section}. Then
$$
\sigma_n\stackrel{n\to\infty}{\lra} \sigma\quad
\text{in}\ \C^0.
$$
\end{proposition}

\proof
We have to show the following. Let $\bar U\times B\to M$ be a local
trivialization with $\bar U\subset B$ a closed ball, and let $V\subset
S^2$ be an open set so that $|C|\cap (\bar U\times V)\to \bar U$ is proper. Let
$d'$ be the degree of $C|_{\bar U\times V}$ over $\bar U$, counted with
multiplicities. Then for $n$ sufficiently large $C_n\cap (\bar U\times
V)\to \bar U$ shall be a (branched) covering of the same degree $d'$. In
fact, any $P\in |C|$ has neighbourhoods of this form with $V$
arbitrarily small. Away from the critical points of the projection to
$\bar U$ both $C$ and $C_n$ would then have exactly $d'$ branches on
$\bar U\times V$, counted with multiplicities. In the coordinates on
$M^{[d]}$ exhibited in Proposition~\ref{M^{[d]}} the components of
$\sigma_n$ are elementary symmetric functions in these branches. As
$V$ can be chosen arbitrarily small this implies $\C^0$-convergence
of $\sigma_n$ towards $\sigma$.

By the definition of the topology on $\shC(M)$ the $C_n$ lie in
arbitrarily small neighbourhoods of $|C|$. Properness of
$|C|\cap(\bar U\times V)\to \bar U$ implies
\begin{eqnarray}\label{boundary inclusion}
\partial(\bar U\times V)\cap |C|\subset \partial \bar U \times V.
\end{eqnarray}
By compactness of $(\partial \bar U\times V)\cap |C|$ we may replace
$|C|$ in this inclusion by a neighbourhood. Therefore (\ref{boundary
inclusion}) holds with $C_n$ replacing $|C|$, for $n\gg0$. We
conclude that $C_n\cap (\bar U\times V)\to \bar U$ is proper for
$n\gg 0$ too, hence a branched covering. Let $d_n$ be its covering
degree.

Convergence of the $C_n$ in the $\C^0$-topology implies that
for every $n$ there exist stable maps $\varphi_n:\Sigma_n\to M$,
$\psi_n:\Sigma_{\infty,n}\to M$ lifting $C_n$ and $C$ respectively
and a $\kappa_n: \Sigma_n\to\Sigma_{\infty,n}$ as above with
$$
d_M\big(\varphi_n(z),\psi_n(\kappa_n(z))\big)\lra 0
$$
uniformly. Let $Z\subset B$ be the union of the critical values of
$p\circ\varphi_n$ and of $p\circ\psi_n$ for all $n$. This is a
countable set, hence has dense complement. Choose $Q\in \bar
U\setminus Z$ and put $F= p^{-1}(Q)$. By hypothesis $F$ is
$J_n$-holomorphic and transverse to $\varphi_n$ and $\psi_n$ for
every $n$. Therefore, for each $n$ the cardinality of
$A_n:=\varphi_n^{-1}(F\cap (\bar U\times V))$ and of
$\psi_n^{-1}(F\cap (\bar U\times V))$ are $d_n$ and $d'$
respectively. Since $P$ is a regular value of $p\circ\psi_n$, for
$n\gg0$ the image $\kappa_n(A_n)$ lies entirely in the regular part
of non-contracted components of $\Sigma_{\infty,n}$. On this part the
pull-back of the Riemannian metric on $M$ allows uniform measurements
of distances. In this metric the distance of $\kappa_n(A_n)$ from
$\psi_n^{-1}(F\cap (\bar U\times V))$, viewed as $0$-cycle, tends to
zero for $n\to \infty$. Therefore $d_n=d'$ for $n\gg0$.
\qed
\medskip

In a situation where the description of
Proposition~\ref{cycles-coefficients},3 applies we obtain convergence
of coefficients, even under the presence of fiber components in the
limit. 

\begin{proposition}\label{coefficient convergence}
Given the data $p:(M,J)\to \CC\PP^1$, $J_0$, $s_0$, $s_1$, $C$, $a_r$
of Proposition~\ref{cycles-coefficients},3 assume that $J_n$ is a sequence
of almost complex structures making $p$ pseudo-holomorphic and so
that $J_n=J_0$ on a neighbourhood of $H\cup p^{-1} \big(p(\bar C\cap
H) \cap p(\bar C\cap S)\big) \cup \bigcup F_a$ that is independent
of $n$. Let $\{C_n\}_n$ be a sequence of $J_n$-holomorphic curves
converging to $C=\bar C+\sum_a m_a F_a$ in the
$\mathscr{C}^0$-topology. Let $a_{0,n}$ be holomorphic sections of
$L^{k+de}$ with zero locus $p(C_n\cap H)$ converging uniformly to
$a_0$.

Then the sections $a_{r,n}$ of $L^{k+(d-r)e}$ corresponding to $C_n$
converge uniformly to $a_r$ for all $r$.
\end{proposition} 

\proof
From Proposition~\ref{cycles-coefficients},3 the $a_{r,n}$ fulfill
equations
$$
\bar\partial a_{r,n}= b_{r,n}(a_{0,n},\ldots,a_{d,n}),
$$
with uniformly bounded right-hand side. Cover $\CC\PP^1$ with 2 disks
intersecting in an annulus $\Omega$ whose closure does not contain
any zeros of $a_0$. Then $H\cap|C| \cap p^{-1}(\Omega)=\emptyset$.
Thus over $\Omega$ the branches of $C_n$ stay uniformly bounded away
from $H$, hence the $a_{r,n}$ are uniformly bounded over $\Omega$.
The Cauchy integral formula on each of the two disks implies a
uniform estimate
$$
\big\|a_{r,n}\big\|_{1,p} \le c\cdot\big(\big\|\bar \partial a_{r,n}\big\|_p
+\big\|a_{r,n}|_\Omega\big\|_\infty\big).
$$
Therefore, in view of boundedness of $b_{r,n}$ everywhere and of
$a_{r,n}$ on $\Omega$ we deduce a uniform estimate on the H\"older
norm
$$
\big\|a_{r,n}\big\|_{0,\alpha}\le c'.
$$
Thus it suffices to prove pointwise convergence of the $a_{r,n}$
on a dense set.

Away from the zeros of $a_0$ union $\bigcup_a p(F_a)$ convergence
follows from Proposition~\ref{section convergence}. In fact, the
quotients $a_{r,n}/a_{0,n}$ occur as coefficients of the local section
$$
s_0^d+p^*(a_{1,n}/a_{0,n})s_0^{d-1} s_1+\ldots+
p^*(a_{d,n}/a_{0,n}) s_1^d=0
$$
of $M^{[d]}\setminus H_d$, see
Proposition~\ref{cycles-coefficients},2. These sections correspond to
a sequence of pseudo-holomorphic curves converging to a
pseudo-holomorphic cycle without fiber components as considered in
Proposition~\ref{section convergence}.
\qed
\medskip

Note that since $\mathscr{C}^0$-convergence $C_n\to C$ implies
convergence of the $0$-cycles $H\cap C_n\to H\cap C$, any sequence
of holomorphic sections $a_{0,n}$ with zero locus $p(C_n\cap H)$
converges after rescaling.

%===========================================================
\section{Unobstructed deformations of pseudo-holomorphic cycles}
\label{sect4}
\noindent
We are interested in finding \emph{unobstructed} deformations of a
pseudo-holo\-mor\-phic cycle $C$ in a pseudo-holomorphic
$S^2$-bundle. In the relevant situations this is possible after
changing the almost complex structure. In this section we give
sufficient conditions for unobstructedness, while the construction of
an appropriate almost complex structure occupies the next section.

The describing PDE follows from Proposition~\ref
{cycles-coefficients},3. Recall the assumptions of loc.cit.: $J$
integrable near $|C|$, standard fiberwise and near $H$ union all
fiber components of $|C|$ union $p^{-1}\big(p(|\bar C|\cap H) \cap
p(|\bar C|\cap S)\big)$. In the notation of loc.cit., to set up the
operator choose $T\subset \O(L^{k+(d-r)e})$, $\mathfrak{D}'\subset
\bigoplus_{r=1}^d L^{k+(d-r)e}$ with
$$
a_0(\CC\PP^1)\times_{\CC\PP^1} \mathfrak{D}'\subset \mathfrak{D}
\quad\text{for all } a_0\in T.
$$
Take $p>2$ and write $W^{1,p}_{\CC\PP^1}(\mathfrak{D}') \subset
\bigoplus_{r=1}^d W^{1,p}(\CC\PP^1, L^{k+(d-r)e})$ for the open set
of Sobolev sections taking values in $\mathfrak{D}'$. View
PDE~(\ref{PDE}) in Proposition~\ref{cycles-coefficients} as family of
differentiable maps
\begin{eqnarray}\label{PDE-map}
\begin{array}{rcl}
W^{1,p}_{\CC\PP^1}(\mathfrak{D}')&\lra& \bigoplus_{r=1}^d
L^p(\CC\PP^1, L^{k+(d-r)e}\otimes\Lambda^{0,1}_{\CC\PP^1}),\\[1ex]
(a_r)_{r=1,\ldots,d}&\longmapsto& \big(\bar\partial a_r-
b_r(a_0,a_1,\ldots,a_d)\big)_{r=1,\ldots,d},
\end{array}
\end{eqnarray}
parametrized by $a_0\in T$. Because the $b_r$ depend holomorphically
on $a_r$ the linearization of this map takes the form
$$
W^{1,p}\big(\textstyle\bigoplus_{r=1}^d L^{k+(d-r)e}\big) \lra
L^p\big(\textstyle\bigoplus_{r=1}^d L^{k+(d-r)e}
\otimes\Lambda^{0,1}_{\CC\PP^1}\big),
\quad v\longmapsto \bar\partial v -B\cdot v.
$$
Here $B$ is a $d\times d$-matrix with entries in $\Hom \big(L^{k+(d-r)e},
L^{k+(d-r')e}\otimes\Lambda^{0,1}_{\CC\PP^1}\big)$.

\begin{proposition}\label{linearization surjective}
Assume that there exists a $J$-holomorphic section $S\subset M$
representing $H-eF$ and that $k\ge 0$. Then $\bar\partial-B$ is
surjective. Moreover, for any $Q\in\CC\PP^1$ the restriction map
$$
\ker\big(\bar\partial -B\big) \lra
{\textstyle \bigoplus_{r\ge1}} L_Q^{k+(d-r)e}
$$
is surjective.
\end{proposition}

\proof
Unlike the case of rank~$1$ surjectivity of $\bar\partial-B$ does not
follow from topological considerations. Instead we are going to
identify the kernel of this operator with sections of the holomorphic
normal sheaf
$$
\shN_{C|M}=\shHom(\shI/\shI^2,\O_C)
$$
of $C$ in $M$, with zeros along $H$. Here $\shI$ is the ideal sheaf of
the possibly non-reduced subspace $C$ of a neighbourhood of $|C|$ in
$M$ where $J$ is integrable. The vanishing of $H^1(\shN_{C|M}(-H))$
then follows essentially topologically. The proof proceeds in three
lemmas.

Local solutions of $\bar\partial v-B\cdot v=0$ form a locally free
$O_{\CC\PP^1}$-module $\shK$ of rank $d$, that is, the sheaf of
holomorphic sections of a holomorphic vector bundle over $\CC\PP^1$ of
rank $d$. Similarly, because $p$ is holomorphic in a neighbourhood of
$C$ the push-forward sheaf
$$
p_*(\shN_{C|M}(-H)): U\longmapsto \shN_{C|M}(-H)
\big(p^{-1}(U)\big)
$$
is a locally free $\O_{\CC\PP^1}$-module of the same rank $d$.

\begin{lemma}\label{shK}
$\shK\simeq p_* \big(\shN_{C|M}(-H)\big)$.
\end{lemma}

\proof
This correspondence holds for any base $B$, so let us write $B$
instead of $\CC\PP^1$ for brevity. Consider first the case $C$ smooth
and transverse to $H$ and $p:C\to B$ having only simple branch points.
Then $\shN_{C|M}$ is the sheaf of holomorphic sections of $\big(
T_M|_C\big) /T_C$. Let $\mathbf{a}$ be the section of $M^{[d]}$
associated to $C$ according to Proposition~\ref{cycles-section}. Away
from the critical points of $p|_C$ the inclusion $T_{M/B}\to T_M$
induces an isomorphism $T_{M/B}|_C\simeq \shN_{C|M}$. Let $\Theta$ be
a holomorphic section of $T_{M/B}(-H)$ along $C\cap p^{-1}(U)$,
$U\subset B$ open. Since $T_{M^d_B/B}=\bigoplus q_i^*(T_{M/B})$,
$q_i:M^d_B\to M$ the $i$-th projection, $\Theta$ induces an
$S_d$-invariant holomorphic section $\tilde\Theta$ of $T_{M^d_B/B}$
over $C\times_B\ldots \times_B C\subset M^d_B$. Recall the
symmetrization map $\Phi: M^d_B\to M^{[d]}$ from Proposition~\ref{cx
structure on M^{[d]}}. By the definition of the almost complex
structure on $M^{[d]}$ this map is holomorphic near $C\times_B\ldots
\times_B C$. Thus $\Phi_*(\tilde\Theta)$ is a holomorphic section of
$\mathbf{a}^*(T_{M^{[d]}/B})$. The vanishing of $\Theta$ on $C\cap H$
translates into the vanishing of the normal component of
$\Phi_*(\tilde\Theta)$ along the divisor $H_d\subset M^{[d]}$
introduced before Proposition~\ref{cycles-coefficients}. In other
words, $\Phi_*(\tilde\Theta)$ is a section of $\mathbf{a}^*\big(
T_{M^{[d]}/B} (-\log H_d) \big)$. It is clear that this sets up an
isomorphism between $p_*(\shN_{C|M}(-H))$ and $\mathbf{a}^*\big(
T_{M^{[d]}/B} (-\log H_d) \big)$. We claim that the module of sections
over $U$ of the latter sheaf is canonically isomorphic to $\shK$. Then
define
$$
\Psi:p_*\big(\shN_{C|M}(-H)\big)\lra \shK
$$
away from the critical points of $p|_C$ by sending $\Theta$ to
$\Phi_*( \tilde\Theta)$.

To prove the claim we have to characterize holomorphic sections of
$T_{M^{[d]}/B}$ along the image of $\mathbf{a}$ in the coordinates
$(z,w): M\setminus (H\cup F)\to \CC^2$. Let $\sigma_1,\ldots,\sigma_d$
be the fiberwise holomorphic coordinates on $M^{[d]}$ induced by $w$,
see Proposition~\ref{cx structure on M^{[d]}}. A section
$\Theta=\sum_r h_r\partial_{\sigma_r}$ of $T_{M^{[d]}/B}$ is
holomorphic iff the $h_r$ depend holomorphically on the $\sigma_i$ and
$[Z,\Theta] =0$, where
$$
Z=\partial_{\bar z}+\beta_1\partial_{\sigma_1}+\ldots+ \beta_d
\partial_{\sigma_d}
$$
is the antiholomorphic horizontal vector field defining the complex
structure. The vanishing of the bracket expresses the fact that $Zf=0$
should imply $Z(\Theta f)=0$. In particular, if $f$ is a holomorphic
function on $M^{[d]}$ then so will be $\Theta f$. Pulling back by
$\mathbf{a}= (\alpha_1(z),\ldots, \alpha_d(z))$ and observing
$\partial_{\bar z}\alpha_r = \beta_r(z,\mathbf{a})$
(Proposition~\ref{cycles-coefficients},2) gives
$$
\partial_{\bar z}\nu_r-{\textstyle \sum_j}
\nu_j\partial_{\sigma_j}\beta_r|_{(z,\mathbf{a})} =0
$$
as condition for $\sum_r \nu_r\partial_{\sigma_r}$ to be a
holomorphic vector field along $\im(\mathbf{a})$. Rescaling the
coordinates $\sigma_r$ by $a_0$ leads to $v=(a_0\nu_1,\ldots,
a_0\nu_d)$. This corresponds to a holomorphic section of
$\mathbf{a}^* \big(T_{M^{[d]}/B} \big)$ iff $\bar\partial v -B\cdot
v=0$.

The map $\Psi$ is an isomorphism wherever defined so far. To finish
the case $C$ smooth, transverse to $H$ and with only simple branch
points it remains to extend over the critical points of $p|_C$. As we
are in a purely holomorphic situation now we are free to work in
actual holomorphic coordinates. Note that if $C$ splits into several
connected components then $M^{[d]}$ is naturally a fibered product
with one factor for each component of $C$. Since our isomorphism
$\Psi$ respects this decomposition the problem is local in $C$. We may
therefore assume $C$ to be defined by $u^2-z=0$ with $u,z$ holomorphic
coordinates and $z$ descending to $B$. In this case $\partial_z$ and
$-u \partial_z$ are a frame for $p_*(\shN_{C|M})$. Away from $z=0$
take linear combinations with $2u\partial_z +\partial_u\in T_C$ to
lift to $T_{M/C}$. Thus
$$
\partial_z=-\frac{1}{2u}\partial_u,\quad
-u\partial_z=\frac{1}{2}\partial_u
$$
in $\shN_{C|M}$. Compute
$$
\begin{array}{rcccccl}
\Psi(\partial_z)&=& -\Psi \big(\frac{1}{2u}\partial_u\big)
&=& -\big( \frac{\sigma_1}{2\sigma_2}\partial_{\sigma_1}
+\frac{\sigma_1^2-2\sigma_2}{2\sigma_2}
\partial_{\sigma_2}\big) \big|_{
\genfrac{}{}{0pt}{}{\sigma_1=0}{\sigma_2=0} }
&=&\partial_{\sigma_2}\\[2ex]
\Psi(-u\partial_z)&=& \Psi \big(\frac{1}{2}\partial_u\big)
&=& \big( \partial_{\sigma_1}
+\frac{\sigma_1}{2} \partial_{\sigma_2}\big) \big|_{
\genfrac{}{}{0pt}{}{\sigma_1=0}{\sigma_2=0} }
&=&\partial_{\sigma_1}.
\end{array}
$$
Therefore $\Psi$ extends to an isomorphism over all of $B$.

For the general case we know by \cite{st99} that locally in $B$ there
exists a holomorphic deformation $\{C_s\}_s$ of $C$ over the
$2$-polydisk $\Delta^2$, say, with $C_s$ smooth and with $p|_{C_s}$
only simply branched for all $s\neq0$. The previous reasoning gives an
isomorphism of holomorphic vector bundles over $B\times
(\Delta^2\setminus \{0\})$. Since the vector bundles extend over
$B\times\{0\}$ this morphism extends uniquely. Define $\Psi$ as the
restriction of this extension to $B\times\{0\}$. Then $\Psi$ is an
isomorphism because the determinant of the extension does not vanish
in codimension $1$; and it is unique because any two deformations of
$C$ fit into a joint deformation with the locus of non-simply branched
curves having higher codimension.
\qed

\begin{lemma}\label{H^1=0}
If $k\ge 0$ then $H^1(C,\shN_{C|M}(-H))=0$ and $p_*(\shN_{C|M}(-H))$
is globally generated.
\end{lemma}

\proof
By Serre duality on $C$
$$
H^1(\shN_{C|M}(-H))^\vee\simeq H^0(C,\shI/\shI^2(H)\otimes K_M(C)) =
H^0(C, K_M(H)|_C).
$$
Now the first part of the statement follows by Lemma~\ref{L has no
sections} below with $L=K_M(H)$. To verify the hypotheses write
$C=\sum_{a=0}^r m_a C_a$ with $C_a=S$ at most for $a=r>0$. If $C_a\sim
d_a H+k_a F$ then for $a<r$
$$
0\le S\cdot C_a=k_a,
$$
and $k_a=0$ only if $d_a>0$. Hence $C_a^2=d_a^2e+2d_a k_a\ge0$ and
$$
\deg_{C_a} K_M(H)= -(d_aH +k_a F)\cdot (H+(2-e)F) = -(2d_a+k_a)<0.
$$
If $C_r=S$ then $\deg_{C_r} K_M(H) =e-2$,
$$
{\textstyle \sum_{a=0}^{r-1}} m_a C_a\sim dH+kF-m_r S
= (d-m_r)H+ (k+m_r e)F,
$$
and for $0\le m<m_r$
$$
m C_r^2 +{\textstyle \sum_{a=0}^{r-1}} m_a C_a\cdot C_r
= -me+(k+m_r e)\ge k+e.
$$
Hence in any case the inequalities required in Lemma~\ref{L has no
sections} hold.

Global generatedness follows by a standard dimension argument with
the Riemann-Roch formula provided $H^1(\shN_{C|M}(-H-F))=0$. This is
also true here because $\deg _{C_a} K_M(H+F)= -(d_a+k_a)<0$ and
$\deg_S K_M (H+F)= e-1$.
\qed
\smallskip

Summarizing, Lemma~\ref{shK} identified the sheaf of local solutions
of $\bar\partial v-B\cdot v=0$ with $p_*(\shN_{C|M}(-H))$. Because
$\bar\partial-B$ is locally surjective standard arguments of
cohomology theory then give an identification
$$
\coker\big(\bar\partial -B\big) \simeq H^1\big(\CC\PP^1,
p_* (\shN_{C|M}(-H))\big).
$$
Because $p|_C$ is a finite morphism the latter sheaf equals $H^1(C,
\shN_{C|M}(-H))$. The latter vanishes by Lemma~\ref{H^1=0}. Moreover,
since by the same lemma $p_*(\shN_{C|M}(-H))$ is globally generated,
so is $\shK$. This gives the claimed surjectivity of the restriction
\qed

\begin{lemma}\label{L has no sections}
Let $C=\sum_{a=0}^r m_aC_a$, $m_a>0$, be a compact holomorphic
$1$-cycle on a complex manifold $X$ of dimension $2$ and let $L$ be a
holomorphic line bundle over $C$. Assume that for all $0<r'\le r$,
$0\le m<m_a$:
$$
c_1(L)\cdot C_0<0,\eqand
c_1(L)\cdot C_{r'}< m C_{r'}^2+{\textstyle \sum_{a=0}^{r'-1}}
m_a C_a\cdot C_{r'}.
$$
Then $H^0(C,L)=0$.
\end{lemma}

\proof
By induction over $\sum m_a$. We identify effective $1$-cycles with
compact complex subspaces without further notice. If $\sum m_a=1$ we
are dealing with a holomorphic line bundle of negative degree over a
reduced space $C=C_0$, which has no non-zero global sections. In the
general case let $\shL$ be the sheaf of holomorphic sections of $L$
and $s\in H^0(C,L)$. The effective cycle $C'=C- C_r$ fulfills the
induction hypothesis. In view of the exact sequence
$$
0\lra \shL\otimes\shI_{C'/C}\lra \shL|_C\lra \shL|_{C'}\lra 0,
$$
the section $s$ lifts to $\shL\otimes\shI_{C'/C}$. Because
$\shI_{C_r}\cdot \shI_{C'}=\shI_C$ the second factor is
$$
\shI_{C'}/\shI_C = \shI_{C'}\otimes \O_X/\shI_{C_r}
= \O_{C_r}(-C').
$$
Thus $\shL\otimes\shI_{C'/C}$ is the sheaf of holomorphic sections of a
line bundle over $C_r$ of degree
\begin{eqnarray*}
\deg_{C_r}\big(\shL\otimes \shI_{C'}/\shI_C\big)
&=& c_1(L)\cdot C_r- C'\cdot C_r\\
&=& c_1(L)\cdot C_r - {\textstyle \sum_{a=0}^{r-1}}
m_a C_a\cdot C_r -(m_r-1)C_r^2,
\end{eqnarray*}
which is $<0$ by assumption. Hence $s$ vanishes identically.
\qed

\begin{remark}
In place of Proposition~\ref{cycles-coefficients} and
Proposition~\ref{linearization surjective} one can use the fact from
complex analytic geometry that the moduli space of compact complex
hypersurfaces in a complex manifold $X$ is smooth at points $C$ with
$H^1(C,\shN_{C|X})=0$. Since this result is not trivial we preferred
to give the elementary if somewhat cumbersome explicit method
described here.
\end{remark}
\smallskip

Surjectivity of $\bar\partial -B$ implies a parametrization of
pseudo-holomorphic deformations of $C$ by a finite-dimensional
manifold. Our unobstructedness result,
Proposition~\ref{unobstructedness of cycle space} below, states this
in a form appropriate for the isotopy problem. In the proof we need
the following version of the Sard-Smale theorem.

\begin{proposition}\label{transverse}
Let $S,X,Y$ be Banach manifolds, $\Phi: S\times X\to Y$ a smooth map
with $\Phi|_{\{s\}\times X}$ Fredholm for all $s\in S$. If $Z\subset
Y$ is a direct submanifold (differential of the inclusion map
splittable) that is transverse to $\Phi$ then the set
\[
\{ s\in S\,|\,\Phi|_{\{s\}\times X}\mbox{\rm\ is transverse to }Z\}
\]
is of second category in $S$.
\end{proposition}

\proof
Apply the Sard-Smale theorem to the projection $\Phi^{-1}(Z)\to S$.
\qed

\begin{proposition}\label{unobstructedness of cycle space}
Let $p:(M,J_0)\to \CC\PP^1$ be a holomorphic $\CC\PP^1$-bundle with
$H,S$ disjoint holomorphic sections, $H\cdot H\ge 0$. For $U,V\subset
M$ open sets with $H\subset V$ consider the space $\mathscr{J}_{U,V}$
of almost complex structures $J$ on $M$ with $J=J_0$ fiberwise and on
$V$, integrable on $U$ and making $S$ holomorphic. Write
$$
\mathscr{M}_{U,V}:=\coprod_{J\in\mathscr{J}_{U,V}} \mathscr{M}_J
$$
with $\mathscr{M}_J$ the space of $J$-holomorphic cycles in $M$.

Assume that $C=\bar C+ \sum m_a F_a$ is a $J$-holomorphic cycle
homologous to $dH+kF$ with $d>0, k\ge 0$, for $J\in
\mathscr{J}_{U,V}$ with $|C|\subset U$, $H\not\subset |C|$, $V$
containing $H\cup p^{-1} \big(p( |\bar C|\cap H)\cap p(|\bar C|\cap
S) \big) \cup \bigcup F_a$. Here $\bar C$ contains all non-fiber
components of $C$. Then
\begin{enumerate}
\item $\mathscr{M}_{U,V}$ and $\mathscr{M}_J$ are Banach
manifolds at $C$.
\item The map $\mathscr{M}_{U,V}\to \mathscr{J}_{U,V}$ is locally
around $C$ a projection.
\item The subset of singular cycles in $\mathscr{M}_J$ is nowhere
dense and does not locally disconnect $\mathscr{M}_J$ at $C$.
Similarly for $\mathscr{M}_{U,V}$.
\end{enumerate}
\end{proposition}

\proof
In Proposition~\ref{linearization surjective} we established surjectivity
of the linearization of the map (\ref{PDE-map}). An application of the
implicit function theorem with $J\in\mathscr{J}_{U,V}$ and $a_0$ as
parameters thus establishes (1) and (2). 

We show the density claim in (3) for $\shM_J$, the case
of $\mathscr{M}_{U,V}$ works analogously. For the time being assume that
$C$ has no fiber components. Apply Proposition~\ref{transverse} with
$S\subset \mathscr{M}_J$ an open neighbourhood of $\{C\}$ and
$$
\Phi: S\times \CC\PP^1\lra M^{[d]},\quad
(C',Q)\longmapsto C'\cap p^{-1}(Q).
$$
For the definition of $M^{[d]}$ see Section~\ref{sect2}. This map
is well-defined for $C'$ close to $C$ by the absence of fiber
components. In local coordinates provided by
Proposition~\ref{cycles-coefficients},1 it is evaluation of
$(a_0,\ldots,a_d)$ at points of $\CC\PP^1$, hence smooth. For
$Z\subset M^{[d]}$ choose the strata $D_i\subset M^{[d]}$ of the
discriminant locus parametrizing $0$-cycles $\sum m_a C_a$ with fixed
partition $d=\sum m_a$ indexed by $i$. The top-dimensional stratum
$D_0$ parametrizes $0$-cycles with exactly one point of multiplicity
$2$; it is a locally closed submanifold of $M^{[d]}$ of codimension
$2$. All other strata $D_i$, $i>0$, have codimension at least $4$.
Since $X$ and $Y$ are finite-dimensional here, the Fredholm condition
is vacuous. Transversality of $\Phi$ to $Z=D_i$ follows for all $i$
from the following lemma.

\begin{lemma}\label{ev is transverse}
For any $Q\in\CC\PP^1$ the map
$$
\mathscr{M}_J\lra S^d(p^{-1}(Q))\simeq \CC\PP^d
$$
is a submersion at $C$.
\end{lemma}
\proof
On solutions of the linearized equation $\bar\partial v -B\cdot v=0$
the differential of the map in question is evaluation at $Q$. The
claim thus follows from the generatedness statement in
Proposition~\ref{linearization surjective}.
\qed
\smallskip

Now since $\codim_\RR D_i>2$ for $i>0$ transversality of
$\Phi|_{\{C'\}\times \CC\PP^1}$ means that $p^{-1}(Q)\cap C'$ has no
point of multiplicity larger than $2$, for all $Q\in\CC\PP^1$. On the
other hand, by  Proposition~\ref{smooth cycles} transversal
intersections with $D_0$ translate into smooth points $P$ of $C'$
with the projection $C'\to\CC\PP^1$ being simply branched at $P$.

For the remaining part of claim (3) we apply
Proposition~\ref{transverse} with $S$ the space of paths
$$
\gamma:[0,1]\to \mathscr{M}_J
$$
connecting two smooth curves $C',C''$ sufficiently close to $C$. The
map is
$$
\Phi: S\times\big([0,1]\times\CC\PP^1)\lra M^{[d]},\quad
(\gamma,t,Q)\longmapsto \gamma(t)\cap p^{-1}(Q).
$$
and $Z$ runs over the $D_i$ as before. Again, transversality follows by
Lem\-ma~\ref{ev is transverse}. For dimension reasons we still obtain
$\gamma(t)\cap D_i=\emptyset$ for $i>0$. It remains to argue
that not only $\gamma$ is transverse to $D_0$ but even $\gamma(t)$ for
every $t\in [0,1]$. Let $W\subset [0,1]\times \CC\PP^1$ be the
one-dimensional submanifold of $(t,Q)$ with $\gamma(t)$ having a point
of multiplicity $2$ over $Q$. Let $v\in T_{\CC\PP^1}$ be in $\ker(Dq)$
where $q$ is the projection
$$
q:W\to [0,1].
$$
Then since $D_0\subset M^{[d]}$ is an analytic divisor and the
differential of $\Phi$ along $t=\mathrm{const}$ is complex linear, it
follows that $i\cdot v$ is also in $\ker(Dq)$. But $W$ is
one-dimensional, so $v=0$ as had to be shown. This finishes the proof,
provided $C$ does not have fiber components.
\smallskip

In the general case let $Z\subset \mathscr{M}_J$ be the subset
of cycles with fiber components. A cycle $C'\in Z$ has a unique
decomposition $C'=\bar C'+\sum m_a F_a$ for certain fibers $F_a$ and
with $\bar C' \sim dH +(k-\sum m_a)F$. The space of such
configurations is of real codimension $2d\sum_a m_a$. Therefore $Z$
is a union of submanifolds of real codimension at least $2$, and
these may be avoided in any path by small perturbations. Apply the
previously established density result to this perturbed path to
obtain a path of smooth cycles.
\qed

%===========================================================
\section{Good almost complex structures}
\label{sect5}
\noindent
Our objective is now to construct an almost complex structure $J$ as
required in Proposition~\ref{unobstructedness of cycle space}, making
an arbitrary pseudo-holomorphic curve in a pseudo-holomorphic
$S^2$-bundle $J$-holomorphic. The next result constructs an
appropriate integrable complex structure $J_0$.

\begin{lemma}\label{integrable complex structure}
Let $p:M\to S^2$ be an $S^2$-bundle and let $J_{M/B}$ be a complex
structure on the fibers of $p$. Let $H,S\subset M$ be disjoint
sections. Then there exists an integrable complex structure $J_0$ on
$M$ with $J_0|_{T_{M/B}}= J_{M/B}$ and making $p$ a holomorphic map
and $H,S$ holomorphic divisors.

Moreover, if $U\subset S^2$ is an open subset and $f: p^{-1}(U)\to
U\times S^2$ is a trivialization mapping $H,S$ to constant sections
then $J_0$ may be chosen to make this trivialization holomorphic.
\end{lemma}

\proof
Since $S\cdot S=- H\cdot H$ we may assume $H\cdot H \ge 0$.
Put $e:=H\cdot H$, and let $F$ be a fiber with $F\subset
p^{-1}(U)$ if $U\neq \emptyset$. It suffices to produce
a map
$$
f:M\setminus (H\cap F)\lra \CC\PP^1
$$
with the following properties.
\begin{enumerate}
\item $f|_{p^{-1}(Q)}$ is a biholomorphism for every $Q\in
S^2\setminus p(F)$.
\item $f^{-1}(0)=H\setminus (H\cap F)$, $f^{-1}(\infty)\subset (S\cup
F)\setminus (H\cap F)$.
\item There exists a complex coordinate $u$ on an open set $U'\subset
S^2$ containing $p(F)$ so that
$$
p^*(u)^e\cdot f: p^{-1}(U')\setminus (H\cup S\cup F)\lra \CC
$$
extends differentiably to a map $p^{-1}(U')\to\CC\PP^1$ inducing a
biholomorphism $F\to \CC\PP^1$.
\end{enumerate}
In fact, away from $F$ this map may be used to define a holomorphic
trivialization of $p$, while near $F$ one may take $u^e\cdot f$.

Since $H\cdot H=e$ a tubular neighbourhood of $H$ is diffeomorphic to
a neighbourhood of the zero section in the complex line bundle of
degree $e$ over $\CC\PP^1$. Let $z,w$ be complex coordinates near
$H\cap F$ with $z=p^*(u)$, $w$ fiberwise holomorphic and $H$ given by
$w=0$. Consider the zero locus of $w- z^e$. This is a local section
of $p$ intersecting $H$ at $P$ of multiplicity $e=H\cdot H$. Hence
this zero locus extends to a section $H'$ isotopic to $H$ and with
$H\cap H'\subset \{P\}$, $H'\cap S=\emptyset$. Under the presence of a
trivialization over $U\subset S^2$ mapping $H,S$ to constant sections
choose $H'$ holomorphic over $U$.

Away from $F$ we now have an $S^2$-bundle with $3$ disjoint sections
and fiberwise complex structure. The uniformization theorem thus
provides a unique map $f:M\setminus F\to \CC\PP^1$ that is fiberwise
biholomorphic and maps $H,H',S$ to $0,1,\infty$ respectively.

It remains to verify (3). Take for $u$ the function with $z=p^*(u)$
as before. Multiplying $f$ by a constant $\lambda\neq 0$ on sections
has the effect of keeping $H$ and $S$ fixed, but scaling $H'=
f^{-1}(1)$ by $\lambda^{-1}$. Thus $u^e\cdot f$ corresponds to the
family with $H'$ replaced by the graph of $z^e/p^*(u^e)\equiv 1$.
This family extends over $u=0$.
\qed
\medskip

We are now in position to construct an almost complex structure so
that a given pseudo-holomorphic cycle has unobstructed deformations.

\begin{lemma}\label{tilde J}
Let $p:(M,J)\to \CC\PP^1$ be a pseudo-holomorphic $S^2$-bundle and
$H,S$ disjoint $J$-holomorphic sections. Let $C\subset M$ be a
$J$-holomorphic curve.

Then for every $\delta>0$ there exist a $\C^1$-diffeomorphism
$\Phi:M\to M$ and an almost complex structure $\tilde J$ on $M$
with the following properties.
\begin{enumerate}
\item $\Phi$ is smooth away from a finite subset $A\subset M$, and
$D\Phi|_A=\id$, $\Phi(S)=S, \Phi(H)=H$.
\item $\|D\Phi\|_{\infty}<\delta$, $\Phi|_{M\setminus
B_\delta(A)} =\id$.
\item $\Phi(C)$ and $H,S$ are $\tilde J$-holomorphic.
\item $p:(M,\tilde J)\to\CC\PP^1$ is a pseudo-holomorphic
$S^2$-bundle.
\item $\tilde J$ is integrable in a neighbourhood of $|C|$.
\item There exists an integrable complex structure $J_0$ on $M$ and
an open set $V\subset M$ containing $H$ and all fiber components of $C$,
so that $\tilde J= J_0$ fiberwise and on~$V$.
\end{enumerate}
\end{lemma}

\proof
Define $F=p^{-1}(\infty)$. Without restriction we may assume
$F\subset C$ and $S,H\not\subset C$. Decompose $C=\bar C\cup
\bigcup_a F_a$ with the second term containing the fiber components.
To avoid discussions of special cases replace $C$ by the closure of
$p^{-1} \big( p(\bar C\cap H) \big)$. For the construction of $J_0$
we would like to apply Lemma~\ref{integrable complex structure}.
However, since we want $\tilde J= J_0$ near the fiber components of
$C$ and $D\Phi|_A=\id$, it is not in general possible to achieve
$J_0|_{ T_{M/\CC\PP^1}}= J|_{T_{M/\CC\PP^1}}$. For each $Q\in
\big(\bar C\cap \bigcup_a F_a\big) \setminus (H\cup S)$ take a local
$J$-holomorphic section $D_i$ of $p$ through $Q$. For each $a$ there
exists a local trivialization $p^{-1}(V_a)= V_a\times\CC\PP^1$
restricting to a biholomorphism $(F_a, J|_{T_{F_a}}) \to\CC\PP^1$ and
sending $D_i$ and $H,S$ to constant sections. Define the fiberwise
complex structure $J_{M/\CC\PP^1}$ near the $F_a$ by pulling back the
complex structure on $\CC\PP^1$ via these trivializations. Extend to
the rest of $M$ arbitrarily. Now define the reference complex
structure $J_0$ by applying Lemma~\ref{integrable complex structure}
with the data $M, J_{M/\CC\PP^1}, S,H$ and the chosen trivialization
near $\bigcup F_a$.

Next we construct the diffeomorphism $\Phi$. Put
$$
A= (C\cup H\cup S)_\sing \cup \mathrm{Crit}(p|_{\bar C_\reg}),
$$
where $\mathrm{Crit}(.)$ denotes the set of critical points of a map.
This is a finite set. For each $P\in A$ let $z,w$ be
$J_0$-holomorphic coordinates near $P$ with $z=p^*(u)$,
$z(P)=w(P)=0$. If $P\in H\cup S$ we also require $w(H \cup S)=\{0\}$.
Let $b(z,w)$ be the function defining $J$ near $P$ via
$T^{0,1}_{M,J}= \langle \partial_{\bar w}, \partial_{\bar z}+ b_P
\partial_w\rangle$. Then
$$
\Psi_P=(z,w- b_P(0,0)\bar w): V_P \lra \CC^2
$$ 
is a chart for $M$ mapping $J|_P$ to the standard complex structure
on $T_{\CC^2,0}$. Note that in this chart $p$ is the projection onto
the first coordinate of $\CC^2$. Now $\Psi_P(C\cap V_P)$ is
pseudo-holomorphic with respect to an almost complex structure
agreeing at $0\in\CC^2$ with the standard complex structure $I$ on
$\CC^2$. Thus Theorem~6.2 of \cite{micallefwhite} applies. It gives a
diffeomorphism $\Phi_P$ of a neighbourhood of the origin in $\CC^2$
of class $\mathscr{C}^1$ mapping $\Psi_P(C\cap V_P)$ to a
holomorphic curve, and with $D\Phi_P|_0 =\id$. Moreover, by our
choices $\Psi_P (H\cup S)$ is already holomorphic and hence, by the
construction in \cite{micallefwhite} remains pointwise fixed under
$\Phi_P$. Therefore there exists, for any sufficiently small $\delta>0$, a
diffeomorphism $\Phi$ of $M$ with
$$
\Phi|_{M\setminus B_\delta(A)}=\id,\quad
\|D\Phi\|<\delta,\quad
\Phi|_{B_{\delta/2}(P)}= \Psi_P^{-1}\circ\Phi_P\circ\Psi_P\quad
\forall P\in A,
$$
and with $\Phi|_H=\id$, $\Phi|_S=\id$. This is the desired
diffeomorphism of $M$.

For the definition of $\tilde J$ observe that on $B_{\delta/2}(P)$,
for $\delta$ sufficiently small and $P\in A$, the transformed curve
$\Phi(C)$ is pseudo-holomorphic with respect to $\Psi_P^*(I)$. On
this part of $M$ define $\tilde J=\Psi_P^*(I)$. Moreover, for $P\in
A\cap (H\cup S\cup \bigcup F_a)$ it holds $b_P\equiv 0$, hence
$\Psi_P^*(I)= J_0$. This is true at $F_a\cap\bar C$ by the definition
of $J_0$, and for $P\in (H\cup S)\cap \bar C$ because $p$ and $H,S$
are pseudo-holomorphic for both $J$ and $J_0$. Therefore we may put
$\tilde J= J_0$ on $B_{\delta/2}(H\cup S\cup \bigcup F_a)$ for
$\delta$ sufficiently small. So far we have defined $\tilde J$ on
$V:= B_{\delta/2}(A\cup H\cup S\cup\bigcup F_a)$. To extend to the
rest of $M$ let
$$
w: M\setminus (H\cup F)\lra \CC
$$
be the restriction of a meromorphic function on $(M,J_0)$ inducing a
biholomorphism on each fiber as in the proof of Lemma~\ref{integrable
complex structure}. Let $u:\CC\PP^1\setminus p(F)\simeq \CC$ and put
$z=p^*(u)$. To define $\tilde J$ agreeing with $J^0$ fiberwise is
equivalent to giving a complex valued function $b$ via
$$
T^{0,1}_{M,\tilde J}=\langle \partial_{\bar w}, \partial_{\bar
z}+b(z,w)\partial_w\rangle,
$$
see Lemmas~\ref{connection and almost cx strs.}, \ref{adapted
coordinates}. The condition that $\Phi(C)$ be pseudo-holomorphic
prescribes $b$ along $\Phi(C_\reg)$. Moreover, $\tilde J$ coincides
with $J^0$ near $H\cup F$ iff $b$ has compact support, and
the already made definition of $\tilde J$ on $V$ forces $b$ to also
vanish there. This fixes $b$ on $\Phi(C)\cup V$. 

\begin{lemma}
There exists an extension of $b$ to $\CC^2$ with compact support and
so that $\partial_{\bar w}b=0$ in a neighbourhood of $\Phi(C)$.
\end{lemma}

\proof
Let $P\in \Phi(C_\reg)$ be a non-critical point of $\Phi(C_\reg)\to
\CC\PP^1$. In a neighbourhood $U_P\subset M$ of $P$ write $\Phi(C)$
as graph $w=\lambda(z)$. We define
$$
b_P(z,w)= \partial_{\bar z} \lambda(z)
$$
on $U_P$. Cover a neighbourhood of $\Phi(C)\setminus V$ with finitely
many such $U_P$. Let $\{\rho_P\}$ be a partition of unity subordinate
to the cover $\{U_P\cap \Phi(C)\}$ of $\Phi(C)\setminus V$. For any $P$ the
projection $p|_{U_P\cap \Phi(C)}$ is an open embedding. Hence there
exists a function $\sigma_P$ on $p(U_P\cap\Phi(C))$ with
$\rho_P=p^*(\sigma_P)|_{\Phi(C)}$. Then $p^*(\sigma_P)|_{U_P}$ is a partition of
unity for $\{U_P\}$ in a neighbourhood $U$ of $\Phi(C)\setminus V$ in
$M$. Put
$$
b(z,w)= \left\{\begin{array}{ll}
\sum_P \sigma_P(z)\cdot b_P(z,w),& (z,w)\in U\\
0&(z,w)\in V.\end{array}\right.
$$
For well-definedness it is crucial that $w$ is globally defined on
$M\setminus (H\cup F)$. Now $b(z,w)$ is a smooth function on a
neighbourhood of $\big(V\cup \Phi(C)\big) \setminus (H\cup F)$ in $M$
with the desired properties. Extend arbitrarily to $M\setminus (H\cup
F)$ with compact support.
\qed
\smallskip

To finish the proof of Lemma~\ref{tilde J} it remains to remark that
$\tilde J$ keeps $p$ pseudo-holomorphic by construction.
\qed
\bigskip

The results of this and the last section will be useful for the
isotopy problem in combination with the following lemma (cf.\ also
\cite{shevchishin}, Lemma~6.2.5).

\begin{lemma}\label{tilde J_n}
In the situation of Lemma~\ref{tilde J} let $\{J_n\}_n$ be a sequence
of almost complex structures making $p$ pseudo-holomorphic and
converging towards $J$ in the $\mathscr{C}^0$-topology on all of $M$
and in the $\mathscr{C}^{0,\alpha}_\loc$-topology on $M\setminus A$.
For every $n$ let $C_n$ be a smooth $J_n$-holomorphic curve with
$C_n\cap A=\emptyset$ and so that $C_n\to C$ in the
$\mathscr{C}^0$-topology. Let $\Phi$, $\tilde J$ be a diffeomorphism
and almost complex structure from the conclusion of Lemma~\ref{tilde
J}.

Then, possibly after going over to a subsequence, there exists a
finite set $\tilde A\subset M$ containing $A$ and almost
complex structures $\tilde J_n$ with the following properties.
\begin{enumerate}
\item $p$ is $\tilde J_n$-holomorphic.
\item $\Phi(C_n)$ is $\tilde J_n$-holomorphic.
\item $\tilde J_n\to \tilde J$ in $\mathscr{C}^0$ on $M$ and in
$\mathscr{C}^{0,\alpha}$ on $M\setminus \tilde A$.
\end{enumerate}

An analogous statement holds for sequences of paths
$\{C_{n,t}\}_t$, $\{J_{n,t}\}_t$ uniformly
converging to $C$ and $J$ respectively.
\end{lemma}

\proof
By the Gromov compactness theorem a subsequence of the $C_n$
converges as stable maps. If $\varphi: \Sigma\to M$ is the limit then
$C=\varphi_* (\Sigma)$. Define $\tilde A$ as the union of $A$ and of
$\Phi\circ\phi$ of the set of critical points of $p\circ \Phi\circ
\varphi$. Note that by the definition of convergence of stable maps,
away from $\tilde A$ the convergence $\Phi(C_n)\to \Phi_*(C)$ is as
tuples of sections. In other words, for $P\in \Phi(|C|)\setminus
\tilde A$, say of multiplicity $m$ in $C$, there exists a
neighbourhood $U_P$ of $P$ and disjoint $J_n$-holomorphic sections
$\lambda_{1,n}, \ldots, \lambda_{m,n}$ of $p$ over $p(U_P)$ so that
$$
\lambda_{i,n}\stackrel{n\to\infty}{\lra} \lambda,\quad i=1,\ldots,m, 
$$
and $\lambda$ has image $\Phi(|C|)\cap U_P$. By elliptic regularity
and $\mathscr{C}^{0, \alpha}$-convergence of the $J_n$ away from $A$
this convergence is even in $\mathscr{C}^{1,\alpha}$. To save on
notation we now write $C$, $C_n$, $J_n$ instead of $\Phi_*(C)$,
$\Phi(C_n)$, $\Phi_*(J_n)$ respectively. The assumptions remain the
same except that $J_n$ may only be continuous at points of $A$.

A diagonal argument reduces the statement to convergence in
$\mathscr{C}^{0,\alpha} (M\setminus B_\eps (\tilde A))$ for any fixed
small $\eps>0$ in place of $\mathscr{C}^{0,\alpha}_\loc$-convergence
on $M\setminus \tilde A$. The construction of $\tilde J_n$ proceeds
on $3$ types of regions, which are $M\setminus B_{\eps/2}(|C|)$,
$B_{3\eps}(\tilde A)$, and $B_\eps(|C|)\setminus B_{2 \eps}(A)$.

On $M\setminus B_{\eps/2}(|C|)$ take $\tilde J_n= \tilde J$.
For the definition near $\tilde A$ observe that $J_n$ fulfills
all requirements except that it is possibly only continuous at
$A$. However, since the distance $\delta_n$ from $C_n$ to $A$ is
positive, there exist smooth $\tilde J_n$ agreeing with $J_n$
away from $B_{\delta_n/2}(A)$ and still converging to $\tilde J$ in
the $\mathscr{C}^0$-topology.

The interesting region is $B_\eps(|C|)\setminus B_{\eps/2}(A)$. Let
$z=p^*(u), w$ be local complex coordinates on $M$ near some
$P\in|C|$ with $u$ holomorphic and $w$ fiberwise $\tilde
J$-holomorphic. Assume that
$$
w=\lambda_{n,i}(z),\quad i=1,\ldots,m, 
$$
describe the branches of $C_n$ near $P$ as above. Since $p$ is
pseudo-holomorphic but $w$ need not be fiberwise pseudo-holomorphic
the almost complex structure $J_n$ is now equivalent to $2$ functions
$a_n, b_n$ via
$$
T^{0,1}_{M,J_n}=\langle \partial_{\bar z}+b_n\partial_w,
\partial_{\bar w}+a_n\partial_w\rangle.
$$
Pseudo-holomorphicity of $\lambda_{n,i}$ means
$$
\partial_{\bar z} \lambda_{n,i} -a_n(z,\lambda_{n,i}) \partial_{\bar
z} \bar\lambda_{n,i} = b_n(z,\lambda_{n,i}).
$$
To define $\tilde J_n$ fiberwise agreeing with $\tilde J$ requires a
function $\tilde b_n$ with
$$
\partial_{\bar z} \lambda_{i,n}= \tilde b_n (z,\lambda_{i,n})
$$
for all $i$.

The intersection of the fibers of $p$ with $\overline{B}_\eps(|C|)$
defines a family of closed disks $\overline{\Delta}_z$ near $P$,
$|z|\ll 1$. Take a triangulation of $\overline{\Delta}_0$ with
vertices $C_n\cap \overline{\Delta}_0$ in the interior of
$\overline{\Delta}_0$ as in the following figure.
\medskip
%-------------------------------------------------------------
\begin{center}
\begin{picture}(0,0)%
\includegraphics{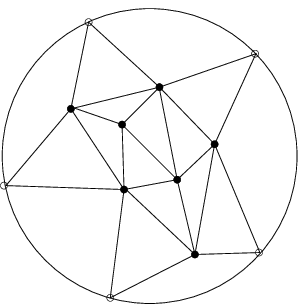}%
\end{picture}%
\setlength{\unitlength}{829sp}%
\begingroup\makeatletter\ifx\SetFigFont\undefined%
\gdef\SetFigFont#1#2#3#4#5{%
  \reset@font\fontsize{#1}{#2pt}%
  \fontfamily{#3}\fontseries{#4}\fontshape{#5}%
  \selectfont}%
\fi\endgroup%
\begin{picture}(6803,6768)(408,-6370)
\put(5536,-2851){\makebox(0,0)[lb]{\smash{\SetFigFont{10}{12.0}{\familydefault}{\mddefault}{\updefault}{\color[rgb]{0,0,0}$\lambda_{i,n}$}%
}}}
\end{picture}
\\
\textsl{\small Triangulation of $C_n\cap \overline{\Delta}_0$.}
\end{center}
\medskip
%-------------------------------------------------------------
The lines in the interior are straight in the $w$-coordinate. Take
$\delta$ so small that the triangulation deforms to
$\overline{\Delta}_z$ for all $|z|\le\delta$. Let $f_{P,n}$ be
the fiberwise piecewise linear function on $B_\eps(|C|) \cap
\{|z|<\delta\}$ restricting to $\partial_{\bar z} \lambda_{i,n} -
\tilde b(z,\lambda_{i,n})$ along the $i$-th branch of $C_n$ and
vanishing at the vertices on the boundary. 

\begin{lemma}
$$
\frac{\big|\nabla (\lambda_{i,n}-\lambda_{j,n})\big|}{
\big|\lambda_{i,n}-\lambda_{j,n}\big|^\alpha} \lra 0.
$$
\end{lemma}

\proof
The difference $u= \lambda_{i,n}-\lambda_{j,n}$ of two branches
fulfills the elliptic equation
\begin{eqnarray*}
\lefteqn{\partial_{\bar z}u -a_n(z,\lambda_{i,n}) \partial_{\bar z}
\bar u}\hspace{1cm}\\
&=& \big( a_n(z,\lambda_{i,n}) - a_n(z,\lambda_{j,n})\big)
\partial_{\bar z} \bar \lambda_{j,n}
+ b_n(z,\lambda_{i,n}) - b_n(z,\lambda_{j,n}).
\end{eqnarray*}
Elliptic regularity gives an estimate
$$
\|\nabla u\|_{0,\alpha} \le c\cdot
\big( \|a_n\|_{0,\alpha} \|\lambda\|_{1,\alpha}+
\|b_n\|_{0,\alpha}+1\big) \cdot \|u\|_{\infty}^\alpha,
$$
with $c$ not depending on $n$. This implies the desired convergence.
\qed
\smallskip

The lemma implies that the H\"older norm of $f_{P,n}$ tends to $0$ for
$n\to \infty$. Let $\tilde f_{P,n}$ be a smoothing of $f_{P,n}$
agreeing with $f_{P,n}$ at the vertices of the triangulation and so
that
$$
\| \tilde f_{P,n}\|_{0,\alpha} \le \|f_{P,n}\|_{0,\alpha} + n^{-1}.
$$
Then near $P$ the desired almost complex structure $\tilde J_n$ will be
defined by $\tilde b_{P,n}= \tilde f_{P,n} +\tilde b$.

To glue keeping $C_n$ and $p$ pseudo-holomorphic we observe that our
local candidates for $\tilde J_n$ fiberwise all agree with $\tilde
J$. It therefore suffices to glue the corresponding sections
$\tilde\beta_n$ of $T^{0,1}_{M/\CC\PP^1}\otimes
p^*(\Lambda^{0,1}_{\CC\PP^1})$ (Lemma~\ref{adapted coordinates})
using a partition of unity. Since Lemma~\ref{adapted coordinates}
requires a fiberwise coordinate $\pi=w$ do this in three steps:
First on a neighbourhood of a general fiber, minus a general section,
then in a neighbourhood of $H$, and finally on $M\setminus (H\cup
F)$.

The statement for paths $\{C_{n,t}\}$, $\{J_{n,t}\}$ follows locally
in $t$ by the same reasoning with an additional parameter $t$;
extend to all $t$ with a partition of unity argument.
\qed

%===========================================================
\section{Generic paths and smoothings}\label{generic paths}\noindent
In this section we discuss the existence of certain generic paths of
almost complex structures. Let $p:(M,\omega,J)\to \CC\PP^1$ be a
symplectic pseudo-holomorphic $S^2$-bundle. For the purpose of this
section $\mathscr{J}$ denotes the space of tamed almost complex
structures on $M$ \emph{making $p$ pseudo-holomorphic}. Endowed with
the $\C^l$-topology $\mathscr{J}$ is a separable Banach manifold.

We will use the following notion of positivity.

\begin{definition}\label{Def. monotone}
An almost complex manifold $(M,J)$ is \emph{monotone} if for every
$J$-holomorphic curve $C\subset M$ it holds
$$
 c_1(M)\cdot C>0.
$$
\end{definition}

\begin{lemma}\label{genericity c_1}
For any $J$ in a path connected Baire subset $\mathscr{J}_\reg\subset
\mathscr{J}$ there exist disjoint $J$-holom\-or\-phic sections
$S,H\subset M$ with $H^2=-S^2\in\{0,1\}$. Moreover, such $J$ enjoy
the following properties.
\begin{enumerate}
\item Any irreducible $J$-holomorphic curve $C\subset M$ not equal to
$S$ is homologous to
$$
d H+ kF,\quad d,k\ge0,
$$
where $F$ is the class of a fiber.
\item $(M,J)$ is monotone (Definition~\ref{Def. monotone}).
\end{enumerate}
In particular, there are no $J$-holomorphic exceptional spheres on $M$
except possibly~$S$.
\end{lemma}

\proof
For $J=I$ the (generic) integrable complex structure $S$ is a
holomorphic section of minimal self-intersection number. Then
$c_1(M,J)$ is Poincar\'e dual to
$$
2S+ (2- S\cdot S) F.
$$
We first consider the case $M=\mathbb{F}_1$. The expected complex
dimension of the space of smooth $J$-holomorphic spheres representing
$S$ is
$$
c_1(M,J)\cdot S+ \dim_\CC(M) -\dim_\CC \Aut(\CC\PP^1) =
(2S+3F)\cdot S-1=0.
$$
This is no surprise as $S$ is an exceptional sphere. For
$J\in\mathscr{J}$ any reducible $J$-holomorphic curve representing $S$
is the union of one section representing $S-lF$ and $l>0$ fibers. The
expected complex dimension of such configurations is
$$
(2S+3F)(S-lF)-1 +l=-l.
$$
Standard transversality arguments show that such reducible
$J$-holomorphic representatives of $S$ do not occur for $J$ in a
path connected Baire subset of $\mathscr{J}$. Note that for any curve
$C\subset M$ variations of $J\in\mathscr{J}$ span the cokernel of the
$\bar\partial_J$-operator, so transversality indeed applies. (If
$C=F$ is a fiber the cokernel is trivial since deformations of $C$ as
fiber span the normal bundle.)

A similar reasoning gives the existence of a $J$-holomorphic section
$H\sim S+F$. Here the expected complex dimension is $2$, so we
impose $2$ incidence conditions to reduce to dimension $0$.

Now if $C\sim dH+kF$ is a $J$-holomorphic curve different from $S$ then
$$
0\le C\cdot S= (dH+kF)\cdot S= k,\quad
0\le C\cdot F= (dH+kF)\cdot F=d
$$
shows that $d,k\ge 0$, $d+k>0$ unless $C=S$. In any case
$$
c_1(M)\cdot C=(2S+3F)\cdot(dH+kF)= 2k+3d>0.
$$

For $M=\CC\PP^1\times\CC\PP^1$ we have $S\cdot S=0$, $c_1(M)=2S+2F$,
the expected complex dimension of spheres representing $S$ is $1$,
and the expected complex dimension of singular configurations
splitting off $l$ fibers is $1-l$. Impose one incidence condition to
reduce the expected complex dimensions by $1$. Proceed as before to
deduce the existence of disjoint $J$-holomorphic sections representing
$S\sim H$ for generic $J$.

This time
$$
0\le C\cdot H=(dH+kF)\cdot H= k,\quad
0\le C\cdot F= (dH+kF)\cdot F=d
$$
shows that $d,k\ge 0$ with at least one inequality strict. Thus again
$$
c_1(M)\cdot C=(2H+2F)\cdot(dH+kF)= 2k+2d>0
$$
for any non-trivial $C\subset M$.
\qed
\bigskip

The other, much deeper genericity result that we will use, is due to
Shevchishin. For the readers convenience we state it here adapted to
our situation, and give a sketch of the proof.

\begin{theorem}\label{genericity Seva}
{\rm (\cite{shevchishin}, Theorems~4.5.1 and 4.5.3.)}
Let $M$ be a symplectic $4$-manifold and $S\subset M$ a finite
subset. There is a Baire subset $\mathscr{J}_\reg$ of the space of
tamed almost complex structures on $M$ with the following properties.
\begin{enumerate} \item $\mathscr{J}_\reg$ is path connected. \item
For a path $\{J_t\}_{t\in[0,1]}$ in $\mathscr{J}_\reg$ let
$\mathscr{M}_{\{J_t\},S}$ be the disjoint union over $t\in[0,1]$
of the moduli spaces of non-multiple pseudo-holomorphic maps
$\varphi:\Sigma\to (M,J_t)$ with $S\subset \varphi(\Sigma)$, for any
closed Riemann surface $\Sigma$. Then there exists a Baire subset of
paths $\{J_t\}$ such that $\mathscr{M}_{\{J_t\},S}$ is a manifold, at
$\varphi$ of dimension
$$
2c_1(M)\cdot \varphi_*[\Sigma] + 2g(\Sigma) -1 - 2\#S,
$$
and the projection
$$
\mathscr{M}_{\{J_t\},S}\lra [0,1]
$$
is open at all $\varphi$ except possibly if $c_1(M)\cdot
\varphi_*[\Sigma]-\sharp S\le 0$, $g(\Sigma)>0$ and $\varphi$ is an
immersion.
\end{enumerate}
\end{theorem}

\proof (Sketch)
We assume $S=\emptyset$, the general case is similar. Let
$\mathscr{M}$ be the disjoint union of the moduli spaces of
$J$-holomorphic maps to $M$ for every $J\in \mathscr{J}$. Define
$\mathscr{M}_\reg$ to be the subset of pairs $(\varphi:\Sigma\to
M,J)$ with cokernel of the linearization of the $\bar
\partial_J$-operator at $\varphi$ having dimension at most $1$. Let
$\mathscr{J}_\reg$ be the complement of the image of
$\mathscr{M}\setminus \mathscr{M}_\reg$ in $\mathscr{J}$. A standard
transversality argument shows that a generic path of almost complex
structures lies entirely in $\mathscr{J}_\reg$. Then one estimates
the codimension of subsets of $\mathscr{M}_\reg$ where the cokernel
is $1$-dimensional. This subset is further stratified according to
the so-called \emph{order} and \emph{secondary cusp index} of the
critical points of $\varphi$. It turns out that a generic path misses
all strata except possibly for the case of order~$2$ and secondary
cusp index~$1$. For this the crucial point is the existence of
explicit second order perturbations of $\varphi$ showing that the map
$$
\nabla D_N: T_{\mathscr{J},J} \lra \Hom(\ker(D_{N,\varphi}),
\coker(D_{N,\varphi}))
$$
is surjective, see \cite{shevchishin}, Lemma~4.4.1. Here
$D_{N,\varphi}$ is a $\bar\partial$-operator on the torsion free part
of $\varphi^*(T_M)/T_\Sigma$. In the immersive case we saw this
operator already in the sketch of the proof of Theorem~\ref{smoothing
nodal curves}. The implicit function theorem then allows to compute
the dimension of strata with $\coker(D_{N,\varphi})$ of specified
dimension. The remaining singularities of order~$2$ and secondary cusp
index~$1$ have local expressions
$$
\varphi(\tau)=\big(\tau^2+O(|\tau|^3), \tau^3+O(|\tau|^{3+\eps})\big),
$$
so are topologically ordinary cusps. Moreover, there must be at least
$c_1(M)\cdot\varphi_*[\Sigma]$ cusps present. 

A further ingredient of the proof is that the presence of a
sufficiently generic cusp contributes complex directions in the
second variation of the $\bar\partial_J$-equation. This implies the
following. Let $\varphi$ be a critical point of the projection
$\mathscr{M}_{\{J_t\}}\lra [0,1]$ and assume
$c_1(M)\cdot\varphi_*[\Sigma]>0$. Then there is a 2-dimensional
submanifold $A\subset\mathscr{M}_{\{J_t\}}$ at $\varphi$ with real
coordinates $x,y$ so that
$$
A\lra \mathscr{M}_{\{J_t\}} \lra [0,1],\quad
(x,y)\longmapsto x^2-y^2.
$$
Therefore $\mathscr{M}_{\{J_t\}} \lra [0,1]$ is open at $\varphi$.
\qed

\begin{remarks}\label{generic J}
1. The proof of the theorem shows that we may restrict to any
subspace $\mathscr{J}$ in the space of tamed almost complex
structures on $M$ having the following properties: For any
$J\in\mathscr{J}$ and any $J$-holomorphic map $\varphi:(\Sigma,
j)\to M$ there exist variations $J_t$ of $J$ in $\mathscr{J}$ and
$j_t$ of $j$ so that terms of the form $(\partial_t J_t,\partial_t
j_t)$ span the cokernel of the linearization of the
$\bar\partial_J$-operator at $\varphi$. Moreover, we need enough
freedom in varying $\partial_t J_t$ in the normal direction to find
solutions of equation~(4.4.6) in \cite{shevchishin}.

Variations of $J$ and $j$ enter in the form $(\partial_t J_t)\circ
D\varphi\circ j+ J\circ D\varphi\circ (\partial_t j_t)$ into this
linearization. Therefore, both conditions are fulfilled if on some
open set of smooth points of $\varphi(\Sigma)$ variations inside
$\mathscr{J'}$ can be prescribed arbitrarily in the normal direction
to $\varphi(\Sigma)$ in $M$.
\smallskip

\noindent
2. Similarly, if variations of $J$ inside $\mathscr{J}$ fulfill the
two conditions only on an open subset $\mathscr{M}'\subset
\mathscr{M}$ then the analogous conclusions of the theorem for
$\mathscr{M}'$ hold true.
\smallskip

\noindent
3. The theorem also holds if we replace $S$ by a path
$\{S(t)\}_{t\in[0,1]}$.
\end{remarks}

\begin{proposition}\label{very generic J}
For $M$ the total space of a symplectic $S^2$-bundle, we may restrict
in the statement of Theorem~\ref{genericity Seva} to almost complex
structures making $p$ pseudo-holomorphic; moreover, we may assume
that for any $J\in\mathscr{J}_\reg$ also the conclusions of
Lemma~\ref{genericity c_1} hold.
\end{proposition}

\proof
By the remark, for $J\in\mathscr{J}$ the only curves that might cause
problems are fibers of $p$. These have always unobstructed
deformations. By the same token, in the definition of
$\mathscr{J}_\reg$ we are free to remove the set of bad almost
complex structures from Lemma~\ref{genericity c_1}.
\qed
\bigskip

The main application of this is the existence of smoothings of
$J$-holo\-mor\-phic cycles occurring along generic paths of almost
complex structures. Our proof uses the unobstructedness of
deformations of nodal curves in monotone manifolds, due to Sikorav.
It generalizes the well-known unobstructedness lemma for smooth
pseudo-holomorphic curves $C$ with $c_1(M)\cdot C>0$ \cite{gromov}
2.1C1, \cite{hoferetal}. In the case where all components are
rational the result is due to Barraud \cite{barraud}. For the readers
convenience we include a sketch of the proof.

\begin{theorem}\label{smoothing nodal curves}
{\rm (\cite{sikorav}, Corollary~2)}
Let $(M,J)$ be an almost complex manifold and $C\subset M$ a
$J$-holomorphic curve with at most nodes and $S\subset C$ a finite
set. Assume that for each irreducible component $C_a\subset C$ it
holds
$$
c_1(M)\cdot C_a> \sharp(C_a\cap S).
$$
Then a neighbourhood of $C$ in the space of $J$-holomorphic cycles is
para\-metrized by an open set in $\CC^d$ with $d=(c_1(M)\cdot C+ C\cdot
C)/2$. The subset parametrizing nodal curves is a union of complex
coordinate hyperplanes. Each such hyperplane parametrizes
deformations of $C$ with one of the nodes unsmoothed.

In particular, for $d>0$ a $J$-holomorphic smoothing of $C$ exists.
\end{theorem}

\proof (Sketch)
We indicate the proof for $\sharp S=\emptyset$, the general case is
similar. Let $\varphi:\Sigma\to M$ be the injective $J$-holomorphic
stable map with image $C$. Standard gluing techniques for
$J$-holomorphic curves give a parametrization of deformations of
$\varphi$ as $J$-holomorphic stable map by the finite-dimensional
solution space of a non-linear equation on $S\times V$; here
$S\subset \CC^N$ parametrizes a certain universal holomorphic
deformation of $\Sigma$ together with some marked points and $V$ is a
linear subspace of finite codimension in a space of sections of
$\varphi^*(T_M)$. The differential of the equation in the
$V$-direction is Fredholm and varies continuously with $s\in S$. The
precise setup differs from approach to approach, see for example
\cite{litian}, \cite{siebert}. The following discussion holds for
either of these.

The statement of the theorem follows from this by two observations.
First, if $\Sigma$ has $r$ nodes then $S$ is naturally a product
$S_1\times S_2$, where $S_1\subset\CC^r$ parametrizes deformations of
the nodes of $\Sigma$, while $S_2$ takes care of changes of the
complex structure of the (normalization of the) irreducible
components of $\Sigma$ together with points marking the position of
the singular points of $\Sigma$. In particular, the $i$-th coordinate
hyperplane in $S_1$ corresponds to deformations of $\Sigma$ with the
$i$-th node unsmoothed. The observation is that the
$\bar\partial_J$-equation is not only differentiable relative $S$ but
even relative $S_1$. In fact, variations along the $S_2$-direction
merely change the complex structure of $\Sigma$ away from the nodes.
This variation is manifestly differentiable in all of the gluing
constructions.

For the second observation let $\hat\Sigma= \coprod_i
\Sigma_i\to\Sigma$ be the normalization of $\Sigma$; this is the
unique generically injective proper holomorphic map with $\hat\Sigma$
smooth. Write $\hat\varphi$ for the composition with $\varphi$. Let
$D_\varphi$ be the linearization of the $\bar\partial_J$-operator
acting on sections of $\varphi^*(T_M)$, and $\widetilde D_\varphi$
the analogous operator with variations in the $S_2$-directions
included. There is a similar operator $D_{N,\varphi}$ acting on
sections of $N:=\hat\varphi^*(T_M)/ d\hat\varphi(T_\Sigma)$, see
e.g.\ \cite{shevchishin}, \S1.5 for details. Both $D_\varphi$ and
$D_{N,\varphi}$ are $0$-order perturbations of the
$\bar\partial$-operator on $\varphi^*(T_M)$ and on $N$ respectively,
for the natural holomorphic structures on these bundles. The adjoint
of $D_{N,\varphi}$ thus has an interpretation as $0$-order
perturbation of the $\bar\partial$-operator acting on sections of
$(N\otimes \Lambda^{0,1}_\Sigma)^*\simeq \det \hat\varphi^*(T^*_M)$.
Thus $\ker(D^*_N)$ restricted to $\Sigma_i$ consists of
pseudo-analytic sections of a holomorphic line bundle of degree
$-c_1(M)\cdot \hat \varphi_*[\Sigma_i]$, which is $<0$ by hypothesis.
On the other hand, every zero of a pseudo-analytic section
contributes positively, see e.g.\ \cite{vekua}. Thus $D_{N,\varphi}$
is surjective.

The point is that this well-known surjectivity implies surjectivity
of $\widetilde D_\varphi$.
Partially descending a similar diagram of
section spaces on $\hat \Sigma$ gives the following commutative
diagram with exact rows:\vfill
\pagebreak

$$
\begin{CD}
0@>>> \mathscr{C}^\infty(T_\Sigma) @>>>
\mathscr{C}^\infty(\varphi^*(T_M)) @>>> \mathscr{C}^\infty(N) @>>>0\\
&& @V{\bar\partial}VV @V D_\varphi VV @VV D_{N,\varphi} V\\
0@>>> \Omega^{0,1}(T_\Sigma) @> d\varphi >>
\Omega^{0,1}(\varphi^*(T_M)) @>>> \Omega^{0,1}(N) @>>>0.
\end{CD}
$$
Here smoothness of a section at a node means smoothness on each
branch plus continuity. Note that the surjectivity of
$\mathscr{C}^\infty(\varphi^*(T_M)) \to \mathscr{C}^\infty (N)$ holds
only in dimension $4$ and because $C$ has at most nodes. The diagram
implies
$$
\coker(D_{N,\varphi})=\Omega^{0,1}(\varphi^*(T_M))\Big/
\Big(d\varphi\big(\Omega^{0,1}(T_\Sigma)\big)+
D_\varphi\big(\mathscr{C}^\infty(\varphi^*(T_M))\big)\Big).
$$
This is the same as $\coker(\widetilde D_\varphi)$. In fact, $\widetilde
D_\varphi$ applied to variations of the complex structure of $\Sigma$
spans $d\varphi(\Omega^{0,1} (T_\Sigma))$ modulo
$D_\varphi(\mathscr{C}^\infty(T_\Sigma))$, and 
$$
D_\varphi(\mathscr{C}^\infty(T_\Sigma)) \subset
D_\varphi(\mathscr{C}^\infty(\varphi^*(T_M))).
$$

Putting together we have shown that $\coker(\widetilde D_\varphi)=0$,
and $\widetilde D_\varphi$ is the differential relative $S_1$ of the
$\bar\partial$-operator acting on sections of $\varphi^*(T_M)$.
\qed
\bigskip

We are now ready for the main result of this section.

\begin{lemma}\label{nodal smoothing}
Let $(M,\omega)$ be a symplectic $4$-manifold and $J$ a tamed
almost complex structure. Assume that $J=J_{t_0}$ for some generic
path $\{J_t\}$ as in Theorem~\ref{genericity Seva},2. Let
$C=\sum_a m_a C_a$ be a $J$-holomorphic cycle with
$c_1(M)\cdot C_a>0$ for all $a$. Suppose that $C_{a'}$ is
an exceptional sphere precisely for $a'>a_0$ and that
$$
C_{a'}\cdot \Big(\sum_{a< a'} m_a C_a\Big) \ge m_{a'}
$$
holds for every $a'>a_0$.

Then there exists a $J$-holomorphic smoothing $C^\dagger$ of $C$.
Moreover, if $U\subset M$ is open and $C|_U= C_1+C_2$ a decomposition
without common irreducible components, then there exists a
deformation of $C$ into a nodal curve $C^\dagger$ with all nodes
contained in $U$ and $C^\dagger\cap U$ the union of separate
smoothings of $C_1$ and $C_2$.
\end{lemma}

\proof
Assume first that no $C_a$ is an exceptional sphere. Let
$\varphi_a:\Sigma_a\to M$ be the generically injective
$J$-holomorphic map with image $C_a$. Consider the moduli space
$\mathscr{M}_{\{J_t\}}$ of tuples of $J_t$-holomorphic maps
$$
(\varphi'_{a,i})_{a,i}:(\Sigma_{1,1},\ldots,\Sigma_{1,m_1},\ldots,
\Sigma_{a,1},\ldots,\Sigma_{a,m_a},\ldots) \lra M^m,
$$
together with $t\in[0,1]$. Here $\Sigma_{a,i}$ denotes the closed
surface underlying $\Sigma_a$ with some complex structure depending
on $i$. By Theorem~\ref{genericity Seva}, $\mathscr{M}_{\{J_t\}}$ is
a smooth manifold of expected dimension and the projection
$$
\mathscr{M}_{\{J_t\}}\lra [0,1]
$$
is open at any $(\varphi'_{a,i})$. The subset of tuples
$(\varphi'_{a,i})$ with one of the maps $\varphi'_{a,i}$ having a
critical point is of real codimension at least $2$. Similarly, since
we have excluded the possibility of exceptional spheres, the
condition that $\varphi'_{a,i}$ and $\varphi'_{a',i'}$ for $(a,i)\neq
(a',i')$ have a common tangent line is also of real codimension $2$.
Thus taking $\{J_t\}$ generic we may assume these configurations to
correspond to a locally finite union of locally closed submanifolds
of real codimension at least $2$. 

We may thus choose a $J$-holomorphic perturbation $(\varphi'_{a,i})$
of the tuple with entries $\varphi_{a,i}=\varphi_a$ avoiding the
subset of singular configurations. The image of $(\varphi'_{a,i})$ is a
$J$-holomorphic curve $C$ with at most nodes and such that each
component evaluates positively on $c_1(M)$. An application of
Theorem~\ref{smoothing nodal curves} now shows that a smoothing of
$C$ exists. This finishes the proof under the absence of exceptional
spheres.

In the general case write $C=\bar C+\sum_e m_e E_e$ with the second
term containing the exceptional spheres. The previous construction
applied to $\bar C$ yields a smooth $J$-holomorphic curve $\Sigma=
\bar C^\dagger\subset M$. Because $\Sigma$ contains no exceptional
sphere and $J=J_{t_0}$, the space of smooth $J$-holomorphic curves
with a common tangent with $E_1$ and homologous to $\Sigma$ is
nowhere dense in the space of all such $J$-holomorphic maps. This
follows by another application of Theorem~\ref{genericity Seva}.
Deforming $\Sigma$ slightly we may therefore assume the intersection
with $E_1$ to be transverse. Now apply Theorem~\ref{smoothing nodal
curves} to $\Sigma\cup E_1$. Because $\Sigma\cdot E_1>0$ the result
is a $J$-holomorphic smoothing $\Sigma^1$ of $\bar C + E_1$ not
containing $E_1$. An induction over $\sum_e m_e$ finishes the
construction. The assumption on the intersection numbers of
exceptional components with the rest of the curve guarantees that in
the induction process the intersection of any exceptional component
with $\Sigma^i$ remains non-empty.

The statement on the existence of partial smoothings is clear from the
proof.
\qed

%===========================================================
\section{Pseudo-holomorphic spheres with prescribed
singularities}\noindent
\begin{proposition}\label{J-hol spheres with prescribed sing's}
Let $p:(M,J)\to \CC\PP^1$ be a pseudo-holomorphic $S^2$-bundle. Let
$W\subset M$ be open and
$$
\varphi:\Delta\lra W
$$
a proper, injective, $J$-holomorphic map. Assume that there exists an
integrable complex structure $J_0$ on $M$ making $p$ a holomorphic
map, and with $J|_W=J_0|_W$, $J|_{T_{M/\CC\PP^1}}=
J_0|_{T_{M/\CC\PP^1}}$.

Then for any $k>0$ there exists a $J$-holomorphic sphere
$$
\psi_k: \CC\PP^1\lra M
$$
approximating $\varphi$ to $k$-th order at $0$:
$$
d_M(\varphi(\tau),\psi_k(\tau)) = o(|\tau|^k).
$$
\end{proposition}
\proof
Let $S,H$ be disjoint holomorphic sections of $p$ with $e:=H\cdot
H>0$, and $F$ a fiber not containing $\varphi(0)$. There exists a
holomorphic map
$$
f: M\setminus (F\cap H) \lra\CC\PP^1
$$
with $f^{-1}(\infty)= H\setminus(F\cap H)$, $f^{-1}(0)\subset S\cup
F$ and inducing an isomorphism on each fiber, cf.\
Lemma~\ref{integrable complex structure}. Let $u:\CC\PP^1\setminus
p(F)\to \CC$ be a holomophic coordinate. Put $z=p^*(u)$ and
$w=f|_{M\setminus (H\cup F)}$.  Thus $z,w$ are holomorphic
coordinates on $M\setminus (H\cup F)$ with $p: (z,w)\mapsto z$ and so
that
$$
T^{0,1}_{M,J}= \langle \partial_{\bar w}, \partial_{\bar z}+
b(z,w)\partial_w\rangle,
$$
see Lemma~\ref{adapted coordinates}. We may assume that the image of
$\varphi$ is not contained in a fiber, for otherwise the proposition
is trivial. By changing the domain of $\varphi$ slightly we may then
assume $(p\circ\varphi)^*(z)=\tau^m$ for the standard coordinate
$\tau$ on the domain $\Delta$ of $\varphi$. 

The Taylor expansion with respect to the coordinates $z,w$ of
$\varphi$ at $\tau=0$ up to order $k$ defines a polynomial map
$\Delta\to M$. For $k\ge m$ this approximation of $\varphi$
takes the form
$$
\tau\longmapsto (\tau^m,h(\tau))
$$
for some polynomial $h$ of degree at most $k$. Since $M$
is a projective algebraic variety, for any $k'>k$ the map
$$
\tau\longmapsto (\tau^m,h(\tau)+\tau^{k'})
$$
extends to a holomorphic map
$$
\psi_0:\CC\PP^1\lra M
$$
osculating to $\varphi$ at $\tau=0$ to order at least $k$. By the
choice of $w$ it holds $\psi_0(\CC\PP^1)\cap S\cap F =\emptyset$.
Hence 
$$
(\psi_0)_*[\CC\PP^1]\cdot S=k',
$$
and $(\psi_0)_*[\CC\PP^1]\sim mH+k' F$. Let $T_{M/\CC\PP^1}$ be the
relative tangent bundle. Because $p$ is pseudo-holomorphic this is a
complex line bundle. Put $l=\deg \psi_0^*(T_{M/\CC\PP^1})$. Now
$c_1(T_{M/\CC\PP^1})$ is Poincar\'e-dual to $2H-eF$, hence
$$
l=(2H-eF)(mH+ k'F)= 2k'+me>k.
$$

To show the existence of $\psi_k$ in the statement of the proposition
we consider a modified Gromov-Witten invariant. Let $\mathscr{J}$ be
the space of almost complex structures $J'$ on $M$ of class $\C^l$
with $J'|_W=J_0|_W$, $J'|_{T_{M/\CC\PP^1}}= J_0|_{T_{M/\CC\PP^1}}$
and making $p$ pseudo-holomorphic. Note that by Lemma~\ref{always
tamed} tamedness on any compact subset in $\mathscr{J}$ is implicit
in this definition. Let $\mathscr{M}$ be the disjoint union over
$J'\in\mathscr{J}$ of the moduli spaces of $J'$-holomorphic maps
$$
\psi:\CC\PP^1\lra M,
$$
with
\begin{enumerate}
\item[(i)] $\psi_*[\CC\PP^1]=(\psi_0)_*[\CC\PP^1]$ in $H_2(M,\ZZ)$.
\item[(ii)] in the chosen holomorphic coordinates $\tau$ on
$\CC\PP^1$ and $(z,w)$ on $M\setminus(H\cup F)$:
$$
\psi(\tau)=(\tau^m,h(\tau)+\tau^{l+1}\cdot v(\tau)).
$$
\end{enumerate}
In particular, $p\circ \psi= p\circ\varphi$.

\begin{lemma}
The forgetful map
$$
\mathscr{M}\lra \mathscr{J}
$$
is a Fredholm map of Banach manifolds of index $0$.
\end{lemma}
\proof
Let $\psi(\tau)=(\tau^m,h+ \tau^{l+1}v_0)$ be $J'$-holomorphic. Let
$b'(z,w)$ be the function describing $J'$ on $M\setminus(H\cup
F)\simeq\CC^2$. Then $\partial_{\bar\tau}= m\bar
\tau^{m-1}\partial_{\bar z}$, and $J'$-holomorphicity of $\psi$ is
equivalent to
$$
\partial_{\bar\tau}(h+\tau^{l+1}v_0)
= m\bar\tau^{m-1} b'(\tau^m, h(\tau)+\tau^{l+1} v_0).
$$
Therefore a deformation $(\tau^m,h+ \tau^{l+1}(v_0+\eta))$ of $\psi$
inside $\mathscr{M}$ is $J'$-hol\-o\-mor\-phic iff
$$
\partial_{\bar \tau}\, \eta= m\frac{\bar\tau^{m-1}}{ \tau^{l+1}}
\Big( b'\big(\tau^m,h+\tau^{l+1}(v_0+\eta)\big) -
b'\big(\tau^m,h(\tau)+\tau^{l+1}v_0\big)\Big). 
$$
Linearizing with fixed $b'$ gives the equation
$$
\partial_{\bar \tau}\, \eta= m\bar\tau^{m-1}(\nabla_w b')\cdot \eta
+ m\bar\tau^{m-1}\Big(\frac{\bar \tau}{\tau}\Big)^{l+1}
(\nabla_{\bar w} b')\cdot \bar \eta.
$$
As $\tau^{l+1}$ globalizes as a section of $\O_{\CC\PP^1}(l+1)$, the
intrinsic meaning of $v$ is as a section of $\psi^* (T_{M/\CC\PP^1})
(-l-1)$. So globally our PDE is an equation of CR-type acting on
sections of a complex line bundle over $\CC\PP^1$ of degree  $-1$.
The index is computed by the Riemann-Roch formula to be $0$.
Moreover, as $\psi$ is generically injective, the cokernel can be
spanned by variations of $b'$. An application of the implicit
function theorem finishes the proof of the lemma.
\qed
\smallskip

Since $\mathscr{J}$ is connected there exists a path $\{J_t\}$ in
$\mathscr{J}$ connecting $J_0$ with $J_1=J$. By a standard
application of the Sard-Smale theorem the restriction of
$\mathscr{M}$ to a \emph{generic} path $\{J_t\}$ is a one dimensional
manifold $\mathscr{M}_{\{J_t\}}$ over $[0,1]$. We claim that
$\mathscr{M}_{\{J_t\}}$ is compact. Let $\psi_i\in\mathscr{M}$ be
$J_{t_i}$-holomorphic, and $t_i\to t_0$. The Gromov compactness
theorem gives a $J_{t_0}$-holomorphic cycle $C_\infty= \sum m_a
C_{\infty,a}$ to which a subsequence of $\psi_i(\CC\PP^1)$ converges.
By the homological condition and the local form of the elements of
$\mathscr{M}$ near $\varphi(0)$ there is exactly one component of
$C_\infty$ that projects onto $\CC\PP^1$. All other components are
fibers. So the expected dimension drops by $2$ for each such bubbling
off. Hence bubbling off can be avoided in a generic path as claimed.

We have thus shown that $\mathscr{M}_{\{J_t\}}\to [0,1]$ is a
cobordism. It remains to observe the following fact.

\begin{lemma}
The fiber of
$$
\mathscr{M}_{\{J_t\}}\lra [0,1]
$$
over $0$ consists of one element.
\end{lemma}

\proof
Let $\psi\in\mathscr{M}_{\{J_t\}}$ be $J_0$-holomorphic. If
$\psi\neq\psi_0$ then the intersection index of $\psi(\CC\PP^1)$ with
$\psi_0(\CC\PP^1)$ at $\psi(0)$ is at least $m\cdot(l+1)$. But
$$
\psi_0(\CC\PP^1)\cdot\psi(\CC\PP^1) = (mH+k'F)^2
=m^2e+2mk'=ml,
$$
which is absurd. Hence $\psi_0$ is the only holomorphic map in
$\mathscr{M}_{\{J_t\}}$.
\qed
\smallskip

We are now ready to finish the proof of Proposition~\ref{J-hol
spheres with prescribed sing's}. The parity of the cardinality of the
fiber stays constant in a one-dimensional cobordism. Hence, by the
lemma, the fiber of $\mathscr{M}_{\{J_t\}}\to [0,1]$ over $1$ must be
non-empty. We may take for $\psi_k$ any element in this fiber. (End
of proof of proposition.)
\qed

%===========================================================
\section{An isotopy lemma}\noindent
Our main technical result runs as follows.

\begin{lemma}\label{isotopy lemma}
Let $p:(M,J)\to \CC\PP^1$ be a pseudo-holo\-mor\-phic $S^2$-bundle
with disjoint $J$-holomorphic sections $H,S$. Let $\{J_n\}$ be a
sequence of almost complex structures making $p$ pseudo-holomorphic.
Suppose that $C_n\subset M$, $n\in\NN$, are smooth $J_n$-holomorphic
curves and that
$$
C_n\stackrel{n\to\infty}{\lra} C_\infty=\sum_a m_a
C_{\infty,a}
$$
in the $\C^0$-topology, with $c_1(M)\cdot C_{\infty,a}>0$ for every
$a$ and $J_n\to J$ in $\C^0$ and in $\C^{0,\alpha}$ away from a
finite set $A\subset M$. We also assume:
\begin{enumerate}
\item[($*$)] If $C'=\sum_a m'_a C'_a$ is a non-zero $J'$-holomorphic
cycle $\mathscr{C}^0$-close to a subcycle of $\sum_{m_a>1} m_a
C_{\infty,a}$, with $J'\in \mathscr{J}_\reg$ as in
Proposition~\ref{very generic J}, then
$$
\sum_{\{a|m'_a>1\}} \big(c_1(M)\cdot C'_a+g(C'_a)-1\big) <
c_1(M)\cdot C'-1.
$$
\end{enumerate}

Then any $J$-holomorphic smoothing $C_\infty^\dagger$ of $C_\infty$
is symplectically isotopic to some $C_n$. The isotopy from $C_n$ to
$C_\infty^\dagger$ can be chosen to stay arbitrarily close to
$C_\infty$ in the $\C^0$-topology, and to be pseudo-holomorphic for a
path of almost complex structures that stays arbitrarily close to $J$
in $\mathscr{C}^0$ and in $\mathscr{C}^{0,\alpha}$ away from a finite
set.
\end{lemma}

\begin{remark}\label{smoothing exists}
A $J$-holomorphic smoothing of $C_\infty$ exists by
Proposition~\ref{nodal smoothing} provided $J=J_{t_0}$ for some
generic path $\{J_t\}_t$ as in Proposition~\ref{genericity Seva}. For the
condition on exceptional spheres observe that for homological reasons
$M$ contains at most one of them, say $S$. Then either $C_\infty=S$
and there is nothing to prove, or $C_n\cdot S\ge 0$ and
$$
S\cdot(C_n- mS)=S\cdot C_n +m\ge m.
$$
\end{remark}

\noindent
\emph{Proof of the lemma.}\ 
By Lemma~\ref{always tamed} there exists a symplectic structure
taming $J$. We may therefore asume $J$ tamed whenever needed.

Recall the functions $m,\delta$ defined in Section~\ref{cycle
topology}. Order the pairs $(m,\delta)$ lexicographically:
$$
(\hat m,\hat \delta)< (m,\delta)\quad\Longleftrightarrow\quad
\hat m< m \quad \text{or}\quad
\big(\hat m= m\eqand \hat \delta<\delta\big).
$$
We do an induction over $(m(C_\infty),\delta(C_\infty))$. The
induction process is inspired by the proof of \cite{shevchishin},
Theorem~6.2.3, but the logic is different under the presence of
multiple components.

The starting point is $(m(C_\infty), \delta(C_\infty))=
(0,0)$. Then $C_\infty$ is a smooth pseudo-holomorphic curve and the
statement is trivial.

So let us assume the theorem has been established for all
$(m,\delta)$ strictly less than $(m(C_\infty), \delta(C_\infty))$.
The proof of the induction step proceeds in 8~paragraphs, referred to
as Steps~1--8.
\medskip

\noindent
1)\ \emph{Enhancing $J$.}\\
One of the two sections $H,S$ deform $J$-holomorphically, say $H$. We
may therefore assume $H\not\subset |C_\infty|$. Apply
Lemma~\ref{tilde J} to $C=|C_\infty| \cup p^{-1}\big(
p(|C_\infty|\cap H) \big)$ and $J$. The result is a
$\C^1$-diffeomorphism $\Phi$ of $M$, smooth away from a finite set
$A'\subset M$, and an almost complex structure $\tilde J$ making $p$
and $H,S$ pseudo-holomorphic. By perturbing $\tilde J$ slightly away
from $|C_\infty|$ we may assume that $(M,\tilde J)$ is monotone.
Moreover, $\Phi_*(C_\infty)$ is now a $\tilde J$-holomorphic cycle
with unobstructed deformation theory in the sense made precise in
Proposition~\ref{unobstructedness of cycle space}. The homological
condition on $C_\infty$ required there follows because $C_n\cdot S\ge
0$. We claim that it suffices to prove the theorem under the
assumption that the limit almost complex structure $J$ has this
special form and \emph{for any particular smoothing
$C_\infty^\dagger$}. By taking the union with the finite set in the
hypothesis we may assume $A'=A$.

Change $C_n$ and $J_n$ slightly to achieve $C_n\cap A =\emptyset$.
For example, apply a diffeomorphism that is a small translation near
$A$ and the identity away from this set. After going over to a
subsequence Lemma~\ref{tilde J_n} now provides almost complex
structures $\tilde J_n$ on $M$ making $p$ pseudo-holomorphic with
$\tilde J_n\to \tilde J$ in $\C^0$ and in $\C^{0,\alpha}$ away from
$A$, and so that $\Phi(C_n)$ is $\tilde J_n$-holomorphic. The
sequence $\Phi(C_n)\to \Phi_*(C_\infty)$ fulfills the hypothesis of
the lemma. Assuming that we can prove the lemma for $J=\tilde J$,
pick for every $n\gg0$ an isotopy $\{\tilde C_{n,t}\}_{t\in[0,1]}$
between $\Phi(C_n)$ and a smoothing of $\Phi_*(C_\infty)$; here
$\tilde C_{n,t}$ is pseudo-holomorphic for a path $\{\tilde
J_{n,t}\}$ with $\tilde J_{n,t}\to \tilde J$ and $\tilde C_{n,t}\to
\Phi_*(C_\infty)$ uniformly in the respective $\C^0$-topologies for
$n\to \infty$. Changing $\tilde C_{n,t}$ and $\tilde J_{n,t}$
slightly near $A$ by an appropriate translation we may also assume
$\tilde C_{n,t}\cap A =\emptyset$ for all $t$. Then
$\{\Phi^{-1}(\tilde C_{n,t})\}$ is an isotopy connecting $C_n$ with a
smoothing of $C_\infty$. Another application of Lemma~\ref{tilde
J_n}, now to $\Phi^{-1}$, allows to find paths $\{J_{n,t}\}$ of
smooth almost complex structures, converging to $J$ uniformly in
$\C^0$ for $n\to \infty$, so that $\Phi^{-1}(\tilde C_{n,t})$ is
$J_{n,t}$-holomorphic.

To get the statement of the lemma, apply this reasoning both to the
original sequence $\{C_n\}$ and to the given sequence of
$J$-holomorphic smoothings of $C_\infty$. This gives an isotopy of
$C_n$ and of some smoothing $C_\infty^\dagger$ with the
$\Phi$-preimage of some $\tilde J$-holomorphic smoothings of
$\Phi_*(C_\infty)$. Finally use the fact that according to
Proposition~\ref{unobstructedness of cycle space}, $\tilde
J$-holomorphic smoothings of $\Phi_*(C_\infty)$ are unique up to
isotopy. This step is crucial to make the connection between $\tilde
J$-holomorphic smoothings and $J$-holomorphic smoothing. Here we need
also Proposition~\ref{coefficient convergence} to the effect that
$\mathscr{C}^0$-convergence of cycles implies convergence of
coefficients in the description of Proposition~\ref
{cycles-coefficients},3.

We can henceforth add the following hypotheses to the lemma.
{\it
\begin{enumerate}
\item[($**$)] For some integrable complex structure $J_0$ on $M$ the
pair $(C,J)$ fulfills the hypotheses of
Proposition~\ref{unobstructedness of cycle space} 
\end{enumerate}}

By Proposition~\ref{unobstructedness of cycle space} it now even
suffices to produce an isotopy of $C_n$ with some $J'$-holomorphic
smoothing for $J'$ sufficiently close to $J$ and of the form as
required in loc.cit. 

In Step~5 we will ask $J$ to have a certain genericity property, which
can be achieved by perturbing $\tilde J$ in the normal direction near
some smooth points of $|C_\infty|$. 
\medskip

\noindent
2)\ \emph{Replacement \, of \, non-multiple \, components \, of \,
$C_\infty$ \, by \, $J$-holomorphic spheres.}\\
Let $P\in |C_\infty|_\sing$ and write
$$
\varphi:\Delta\lra M
$$
for any pseudo-holomorphic parametrization of a branch of
$|C_\infty|$ that belongs to a non-multiple component of $C_\infty$
at $P$. Denote by $\bar C_\infty$ the sum of the multiple components
of $C_\infty$, as cycle. For the singularities at $P$ of the curves
underlying $C_\infty$ and $\bar C_\infty$ we have
$$
(|C_\infty|,P)= (|\bar C_\infty|,P) \cup
{\textstyle \bigcup_\varphi}(\varphi(\Delta),P).
$$
Since $J=J_0$ near $P$, by hypo\-thesis~($**$) from Step~1,
Proposition~\ref{J-hol spheres with prescribed sing's} applies. We
obtain a sequence
$$
\psi_k:\CC\PP^1\lra M
$$
of $J$-holomorphic spheres approximating $\varphi$ at $0\in
\Delta\subset \CC\PP^1$ to order $k$. Note that the approximation of
fiber components is exact. We do not indicate the dependence of
$\psi_k$ on $\varphi$ in the notation. For $k\gg0$ the topological
types of the curve singularities $(|C_\infty|,P)$ and of $(|\bar
C_\infty|,P) \cup \bigcup_\varphi (\psi_k(\CC\PP^1),P)$ agree. Fix
such a $k$. For $\delta>0$ define the annulus
$A_\delta=B_{2\delta}(P)\setminus B_{\delta}(P)$. For $\delta>0$
sufficiently small the intersections
$$
A_\delta\cap |C_\infty|\eqand A_\delta\cap (|\bar C_\infty|\cup
{\textstyle \bigcup_{\varphi}} \psi_k(\CC\PP^1))
$$
are isotopic by an isotopy in $A_\delta$ that is the identity near
$|\bar C_\infty|$. Therefore, for each non-multiple branch $\varphi$
there exists a map
$$
\psi:\CC\PP^1\lra M,
$$
which agrees with $\psi_k$ over $M\setminus B_{2\delta}(P)$ and with
$\varphi$ over $B_\delta(P)$, and which is isotopic to $\varphi$ on
$A_\delta$. Let $J'$ be a small perturbation of $J$ agreeing with $J$
near $|\bar C_\infty|$ and making $p$ and $\psi$ pseudo-holomorphic. 
At the expense of further changes of $J$ and $\psi$ away from $|\bar
C_\infty|$, we obtain an almost complex structure $J'$ and a
$J'$-holomorphic map $\psi':\CC\PP^1 \to M$ with
\begin{itemize}
\item $J'$ agrees with $J$ in a neighbourhood of $|\bar C_\infty|$.
\item $\psi'$ is isomorphic to $\varphi$ over a neighbourhood of
$P$.
\item Except possibly at $P$ the map $\psi'$ is an immersion intersecting
$|\bar C_\infty|$ transversely. 
\end{itemize}

Do this inductively for all reduced branches of $C_\infty$ at all
singular points of $|C_\infty|$. The result is an almost complex
structure $J'$ agreeing with $J$ near $|\bar C_\infty|$ and making
$p$ pseudo-holomorphic, and a $J'$-holomorphic cycle $C_\infty'$ all
of whose non-multiple components are rational. We refer to the part
of the added components away from the singularities of $|C_\infty|$
as \emph{parasitic}. Let $U\subset M$ be a neighbourhood of $|\bar
C_\infty|$ so that $C_\infty'$ agrees with $C_\infty$ on $U$ except
for smooth branches of parasitic components \emph{away} from the
singularities of $|C_\infty|$. We can take the intersection of $U$
with the parasitic part of $|C'_\infty|$ to be a union of disks.

Similarly adjust the sequence $C_n$. By $C^{0,\alpha}$-convergence of
the almost complex structures the tangent spaces of $C_n$ converge to
the tangent spaces of $|C_\infty|$ away from the multiple components
of $C_\infty$. Hence, for all $n\gg0$ there exist an almost complex
structure $J'_n$ making $p$ pseudo-holomorphic, a $J'_n$-holomorphic
immersion $\varphi'_n:\Sigma'_n \to M$ and an open set $V_n\subset
\Sigma'_n$ with the following properties.
\begin{itemize}
\item $J'_n=J_n$ on $U$, $J'_n
\stackrel{n\to\infty}{\lra} J'$ in $\C^0$ and in
$\C^{0,\alpha}$ away from $|C_\infty|_\sing\cup A$.
\item $\varphi'_n(\Sigma'_n) \stackrel{n\to\infty}{\lra} C'_\infty$.
\item $\varphi'_n|_{V_n}$ is isomorphic to the inclusion $C_n\cap
U\to M$, which is a part of the original curve, as
pseudo-holo\-mor\-phic map.
\end{itemize}
Define $C'_n=\varphi'_n(\Sigma'_n)$. We call $\Sigma'_n\setminus
V_n$ and $\varphi'_n(\Sigma'_n\setminus V_n)$ the \emph{parasitic}
part of $\Sigma'_n$ and of $C'_n$ respectively.

A note on notation: Because for the most part of the proof we work
with curves having parasitic part we henceforth drop the primes on
all symbols. To refer to the original curves and almost complex
structures we place an upper index $0$, so the original sequences now
read $C_n^0 = \varphi_n^0(\Sigma_n^0)\to C_\infty^0$ and $J_n^0\to
J^0$.
\medskip

\noindent
3)\ \emph{Description of further strategy.}\\
If $C_\infty$ has multiple components the further strategy is to
deform $\varphi_n$ as pseudo-holomorphic map \emph{away} from
$C_\infty$, uniformly in $n$. In the deformation process we fix
enough points so that any occurring degeneration has better
singularities, measured by $(m,\delta)$. Thus induction applies and we
can continue with any smoothing of the degeneration. Therefore, the
deformation process is always successful. The result is a sequence
of $J_n$-holomorphic curves $\{C'_n\}$ with $C'_n$ containing the set
of chosen points and symplectically isotopic to $C_n$ in the required
way, but with a uniform distance from $|C_\infty|$. By the same
reason as before $\lim C'_n$ is a $J$-holomorphic cycle with better
singularities. Next we replace the parasitic part with the removed
part of $C^0_\infty$ and apply the induction hypothesis, if
necessary. This will eventually lead to a $J'$-holomorphic smoothing
$(C^0_\infty)^\dagger$ of the original cycle $C^0_\infty$ that is
isotopic to $C^0_n$, with $J'$ arbitrarily close to $J$ and
fulfilling the other conditions requested in Proposition~\ref
{unobstructedness of cycle space}. This was left to be shown in
Step~1.

If $C_\infty$ is reduced this argument does not work because we can
not prescribe enough points in the deformation to move away. Instead
the pseudo-holomorphic deformation $\varphi_{n,t}$ of $\varphi_n$ will
now be with respect to a generic path $\{J_{n,t}\}_t$ of almost complex
structures connecting $J_n$ with $J$. The deformation
process is successful as long as $\varphi_{n,t}$ stays close to
$\varphi_n$. Otherwise a diagonal argument gives a sequence
$\varphi_{n,t_n}$ converging to a $J$-holomorphic curve with improved
singularities as above.

Alternatively, the reduced case follows from the local isotopy
theorem due to Shevchishin \cite{shevchishin}, Theorem~6.2.3. In
fact, in this case our proof comes down to a global version of the
proof given there. For the sake of completeness we nevertheless
provide full details here.
\medskip

\noindent
4)\ \emph{Replacements of parasitic parts.}\\
The rest of the proof will repeatedly require replacements of the
parasitic part by the removed part of the original curves $C_n^0$ or
$C_\infty^0$. The purpose of this paragraph is to make this process
precise. A point of caution concerns the nodes that have been
generated by the introduction of parasitic components. To simplify
the proof we use the fibered structure here, although this is not
strictly necessary.

Denote by $C_\infty^{\mathrm{par}}$ the union of the non-multiple
irreducible components of $C_\infty$ \emph{that are not fibers of
$p$}. Away from $p(|C_\infty|_\sing)$ the projection
$$
p: C_\infty^{\mathrm{par}}\lra \CC\PP^1
$$
is a holomorphic (branched) covering. For small $\delta'$
$$
A:= B_{2\delta'}\big(p(|C_\infty^0|_\sing)\big)\setminus
B_{\delta'}\big(p(|C_\infty^0|_\sing)\big)
$$
is a union of annuli containing no branch points of this covering.
Take $\delta'$ also so small that $C_\infty^{\mathrm{par}} \cap
p^{-1}(A)\Subset U$, $U$ the neighbourhood of $\bar C_\infty =\bar
C_\infty^0$ from Step~2. Then
$$
\tilde A:=C_\infty^{\mathrm{par}} \cap p^{-1}(A)\cap U
$$
is also a union of annuli, one for each branch of
$C_\infty^{\mathrm{par}}$ near a singular point of $|C^0_\infty|$. Choose
$\delta''< d(\tilde A, M\setminus U)$ and less than
$$
\frac{1}{2} \min \big\{ d(x,y)\big| x,y\in p^{-1}(A)\cap
|C_\infty^0|\cap U, x\neq y, p(x)=p(y) \big\}.
$$
Then there exists a tubular neighbourhood
$$
\Theta: \tilde A\times \Delta\lra M,
$$
of $\tilde A$ with image a union of connected components of
$B_{\delta''} (C_\infty^{\mathrm{par}})\cap p^{-1}(A)$ and so that
$\Delta_Q:=\Theta(\{Q\} \times \Delta)$, $Q\in \tilde A$, are
$J$-holomorphic disks contained in the fibers of $p$. By construction
it holds $\Delta_Q\cdot C_\infty=1$ for the intersection number,
hence also $\Delta_Q\cdot C_n=1$ for $n\gg0$. Therefore,
$\Theta^{-1}(C_n)$ is the graph of a function $\tilde A\to \Delta$.
The same reasoning applies to any small pseudo-holomorphic
deformation of $\varphi_n$ with respect to almost complex structures
keeping $p$ pseudo-holomorphic. Recall from Step~2 that the gluing
took place outside of $U$. We will reglue on $\im(\Theta)$, that is,
away from a smaller neighbourhood $U'\subset U$ chosen disjoint
from $\im(\Theta)$ now. Denote by $d_\cyc$ a metric on $\operatorname
{Cyc}_{\mathrm{pshol}}(M)$ inducing the $\mathscr{C}^0$-topology.
Normalize $d_{\cyc}$ in such a way that $d_\cyc(\hat C,C_\infty)\le
\eps$ implies $|\hat C|\subset B_\eps(|C_\infty|)$.

\begin{lemma}\label{replace parasitic}
For $\eps\ll1$ and $n\gg0$ the following holds. Let
$\{J_{n,t}\}_{0\le t\le t_0}$, be a continuous family of almost
complex structures with $p:(M,J_{n,t})\to \CC\PP^1$
pseudo-holomorphic, $J_{n,0}=J_n$ and
$\|J_{n,t}-J_{n,t'}\|_{0,\alpha}<\eps$ for all $t,t'$. Assume that 
$\big\{\varphi_{n,t}: (\Sigma_n,j_t)\to M\big\}_t$ is a continuous
family of $J_{n,t}$-holomorphic immersions with
$\varphi_{n,0}=\varphi_n$, injective over $B_\delta
(|C^0_\infty|_\sing)$ and such that
$$
d_\cyc\big(\varphi_{n,t}(\Sigma_n),C_\infty\big)\le \eps.
$$

Then there exist a continuous family $\{J^0_{n,t}\}_t$ of almost
complex structures making $p$ pseudo-holomorphic and a continuous
family of smooth $J^0_{n,t}$-holomorphic curves $\{C^0_{n,t}\}_t$
with the following properties.
\begin{enumerate}
\item $C^0_{n,0}=C^0_n$, $J_{n,0}^0=J_n^0$.
\item $d_\cyc(C^0_{n,t},C^0_\infty)\le 2\eps$,
$|\!|J^0_{n,t}- J^0_{n,t'}|\!|_{0,\alpha} < 2\eps$\, for all $t,t'$.
\item $C^0_{n,t}\cap (M\setminus U) = C^0_n\cap (M\setminus U)$.
\item $C^0_{n,t}\cap U'$ is a union of irreducible components of
$\im(\varphi_{n,t})\cap U'$ containing $\im(\varphi_{n,t})\cap U' \cap
B_\delta(|C_\infty^0|_\sing)$.
\end{enumerate}
\end{lemma}

\proof
For $\eps\ll1$ and $n\gg 0$ the preimage of $\partial(\im(\Theta))$
under $\varphi_{n,t}$ is a union of smooth circles
$\Lambda_t\subset\Sigma_n$, moving continuously with $t$. Since
$\im(\Theta)\subset U$, the parasitic part of $\varphi_n$
lies in the complement of $\Lambda_0$. By the construction of the
parasitic part any connected component of $\Sigma_n\setminus
\Lambda_0$ containing some parasitic part is a disk. At time $t$
remove all connected components from $\Sigma_n \setminus\Lambda_t$
that are deformations of such a parasitic component. Call the
resulting map $\bar\varphi_{n,t}$. The assumptions imply that
$\bar\varphi_{n,t}$ is an embedding, provided $\eps>0$ is sufficiently
small and $n$ sufficiently large.

Now for $\eps\ll1$ and $n\gg 0$, on $\im(\Theta)$ the $J_{n,t}$ get
arbitrarily close to $J$ in the $\C^{0,\alpha}$-topology, and the
image of $\varphi_{n,t}$ is a graph over $\tilde A$. Elliptic
regularity thus allows to bound the slope of $\Theta^{-1}\big
(\im(\varphi_{n,t}) \big)$ in terms of $d_\cyc(C_{n,t}), C_\infty)$.
Hence we can glue $\im(\bar\varphi_{n,t})$ with the reduced part of
$C^0_n$ at only a slight change of slopes to obtain $C^0_{n,t}$. 
Therefore, for $\eps\ll1$ and $n\gg0$ there exist $J^0_{n,t}$,
$C^0_{n,t}$ enjoying the requested properties.
\qed
\medskip

\noindent
5)\ \emph{The incidence conditions $ |\mathbf{x}|\subset
|C_{\infty}|$.}\\
Write
$$
C_\infty=\sum_a m_a C_{\infty,a}.
$$
(Note that by our convention from the end of Step~2 the $m_a$ and
$C_{\infty,a}$ in the hypothesis of the lemma are now $m^0_a,
C^0_{\infty,a}$.) For each $a$ choose points
$x^a_1,\ldots,x^a_{k_a}\in C_{\infty,a}$ in general position with
$$
k_a:=c_1(M)\cdot C_{\infty,a}+g(C_{\infty,a})-1.
$$
Since $c_1(M)\cdot C_{\infty,a}>0$ this number is non-negative. Note
that hypothesis~($*$) applied to $C'=\bar C_\infty$ together with
rationality of the reduced part of $C_\infty$ implies
\begin{eqnarray}\label{number of constraints}
k:=\sum_a k_a\le c_1(M)\cdot C_\infty-1.
\end{eqnarray}
The inequality is strict if $C_\infty$ has multiple components. In
the reduced case all components of $C_\infty$ are rational and then
equality holds. Write $\mathbf{x}^a= (x^a_1,\ldots,x^a_{k_a})$ and
$\mathbf{x}=(x_1,\ldots,x_k)$ for the concatenation of the
$\mathbf{x}^a$. Let $\varphi_{\infty,a}: \Sigma_{\infty,a}\to M$ be
the generically injective $J$-holomorphic map with image
$C_{\infty,a}$. Consider the moduli space
$\mathscr{M}_{a,\mathbf{x}^a}$ of pseudo-holomorphic maps
$\varphi:(\Sigma_{\infty,a},j_a)\to M$ representing $[C_{\infty,a}]$
in homology and with $|\mathbf{x}^a|\subset \varphi
(\Sigma_{\infty,a})$. Here $j_a$ may vary. By Remark~\ref{generic
J},2 we may assume $J$ to be regular for all small, non-trivial
deformations of $\varphi_{\infty,a}$ keeping the incidence conditions
$\mathbf{x}^a$; it suffices to change the almost complex structure in
the normal direction near a general point of $C_{\infty,a}$. If
$m_a>1$ this leads to a change of $\tilde J$ in Step~1; if $m_a=1$
change $J'$ away from the neighbourhood $U$ of $|\bar C_\infty|$ in
its construction in Step~2.

Now the dimension of $\mathscr{M}_{a,\mathbf{x}^a}$ is
$$
\big(c_1(M)\cdot C_{\infty,a}+g(\Sigma_{\infty,a})-1\big)-k_a=0.
$$
Hence $\varphi_{\infty,a}$ is an isolated point of
$\mathscr{M}_{a,\mathbf{x}^a}$ for every $a$. This implies that the
singularities of $C_\infty$ under any non-trivial deformation as
$J$-holomorphic cycle containing $\mathbf{x}$ improve. We measure
this improvement through the pair $(m,\delta)$. Let
$\delta^{\mathrm{par}}$ be the number of virtual double points of
$|C_\infty|$ restricted to a neighbourhood of the parasitic part. Because
the parasitic part intersects $|\bar C^0_\infty|$ transversally
it holds
$$
(m(C_\infty),\delta(C_\infty)) =(m(C^0_\infty),
\delta(C^0_\infty)+\delta^{\mathrm{par}}).
$$

\begin{proposition}\label{no equisingular deformations of C_infty}
Let $\hat C$ be a $\hat J$-holomorphic cycle sufficiently close to
$C_\infty$ with $|\mathbf{x}|\subset |\hat C|$.
Assume that one of the following applies.
\begin{enumerate}
\item[(i)] $\hat J= J_t$ for a general path $\{J_t\}$ of almost complex
structures and $|\hat C|$ is not a nodal curve.
\item[(ii)] $\hat J=J$, $\hat C\neq C_\infty$.
\end{enumerate}
Then
$$
(m(\hat C),\delta(\hat C))<(m(C_\infty),\delta(C_\infty)).
$$
\end{proposition}

\proof
Lemma~\ref{semicontinuity of m,delta} implies
$(m(\hat C),\delta(\hat C))\le (m(C_\infty),\delta(C_\infty))$
for $\hat C$ sufficiently close to $C_\infty$ in the $\C^0$-topology.
Moreover, if equality holds for some $\hat C$ containing
$|\mathbf{x}|$ then every component of $\hat C$ is the image of some
$\hat J$-holomorphic deformation $\varphi_a$ of $\varphi_{\infty,a}$
with $|\mathbf{x}^a| \subset \im(\varphi_a)$. Inside the
para\-me\-triz\-ing moduli space of generically injective
pseudo-holomorphic maps with incidence conditions the subspace of
non-immersions is of real codimension at least $2$. Hence such maps do
not occur over a general path of almost complex structures.

In the second case $\hat J=J$ we saw before that $\varphi_{\infty,a}$
is an isolated point in $\mathscr{M}_{a,\mathbf{x}^a}$. So equality
can only occur if $\varphi_a = \varphi_{\infty,a}$ for all $a$, that
is $\hat C=C_\infty$.
\qed
\medskip

\noindent
6)\ \emph{The non-reduced case: Deforming $C_n$ uniformly away from
$C_\infty$.}\\
Deform $C_n$ and $J_n$ slightly to achieve $|\mathbf{x}| \subset C_n$
for all $n$. For example, near $x_i$ apply a locally supported
diffeomorphism of $M$ moving one branch of $C_n$ to $x_i$. As before
let $\varphi_n:(\Sigma_n,j_n) \to M$ be the generically injective
$J_n$-holomorphic map with image $C_n$.

We first treat the case when $C_\infty$ has multiple components.
The purpose of this section is to establish the following result.

\begin{lemma}\label{C'_n}
For every sufficiently small $\eps>0$ there exists, for every $n\gg0$, a
$J_n$-holomorphic immersion $\varphi'_n:(\Sigma_n,j'_n)\to M$ with
\begin{enumerate}
\item $\varphi'_n$ is homotopic to $\varphi_n$ through a family of
immersions, pseudo-holo\-mor\-phic with respect to a path of almost
complex structures staying arbitrarily close to $J_n$ in
$\mathscr{C}^{0,\alpha}$.
\item $|\mathbf{x}|\subset \varphi'_n(\Sigma_n)$.
\item $d_\cyc(\im(\varphi'_n),C_\infty)=\eps$.
\end{enumerate} 
\end{lemma}

\proof
The proof is by descending induction over the degree $c_1(M)\cdot
C_n$. For a first try take paths $x_n(t)$ on $M$ with $x_n(0)\in
C_n\setminus |\mathbf{x}|$ and $\frac{d}{dt} d_M(x_n(t),|C_\infty|)
=2\eps$. Such paths exist for all $n$ provided that $\eps$ is
sufficiently small. Write $\mathbf{x}_n(t)$ for the concatenation
$(x_1,\ldots, x_k,x_n(t))$ of $\mathbf{x}$ and $x_n(t)$. Choose a
$\mathscr{C}^{0,\alpha}$-small deformation $\{J_{n,t}\}_{t\in[0,1]}$
of $J_n$ so that
\begin{itemize}
\item $\{J_{n,t}\}$ fulfills the properties of Theorem~\ref{genericity
Seva},2 with $S(t)=|\mathbf{x}_n(t)|$.
\end{itemize}
Note that the expected (real) dimension of deformations of $C_\infty$
not decreasing $(m,\delta)$ and containing $|\mathbf{x}_n (t)|$ is
$-2$, so these do not occur for generic $\{J_{n,t}\}$, cf.\
Proposition~\ref{no equisingular deformations of C_infty}.

Now try to deform $\varphi_n$ as $J_{n,t}$-holomorphic map with
$|\mathbf{x}_n(t)|$ contained in the image. By
Theorem~\ref{genericity Seva} the set of $t\in[0,1]$ where this is
possible, is open. For closedness let $C_{n,t_i}$ be the image of a
deformation of $\varphi_n$ at times $t_i\to t_\infty$. By the Gromov
compactness theorem we may assume cycle-theoretic convergence
$$
C_{n,t_i}\lra C_{n,t_\infty}.
$$
If $C_{n,t_\infty}$ is reduced and irreducible Lemma~\ref{nodal
smoothing} gives a $J_{n,t_\infty}$-holomorphic partial smoothing of
$C_{n,t_\infty}$ containing $|\mathbf{x}_n(t)|$, see also
Remark~\ref{smoothing exists}. The smoothing has to be partial
because we want to keep deformations of the parasitic components
$C_\infty^{\mathrm{par}} \cap (U'\setminus B_{\delta}
(|C_\infty^0|)$. It remains to verify that this partial smoothing is
symplectically isotopic to $C_{n,t_i}$ for $i\gg0$. Replace the
parasitic part by the original curve, both for the sequence
$C_{n,t_i}$ and for the partial smoothing. This is possible by
Lemma~\ref{replace parasitic}. Since $C_{n,t_\infty}$ is reduced the
induction hypothesis provides a symplectic isotopy between the curves
with parasitic parts replaced. Glue back the parasitic parts
throughout the isotopy to obtain the requested homotopy of
immersions.

If $C_{n,t_\infty}$ is non-reduced or reducible write
$$
C_{n,t_\infty}=\sum_{a=1}^{a_0} C_{n,t_\infty,a}.
$$
Here the $C_{n,t_\infty,a}$ may not be distinct. Split
$\mathbf{x}$ arbitrarily into $a_0$ tuples
$|\mathbf{x}^a|$ with
$$
|\mathbf{x}^a|\subset C_{n,t_\infty,a}.
$$
Let $k_a$ be the number of entries of $\mathbf{x}^a$. Since by
(\ref{number of constraints}) the number of entries $k=\sum_a k_a$ of
$\mathbf{x}$ is less than $c_1(M)\cdot C_{n,t_\infty}-1$ there
exists an $a$ with
$$
k_a< c_1(M)\cdot C_{n,t_\infty,a}-1.
$$
Let $\varphi'_{n,t_\infty}$ be the generically injective map with
image $C_{n,t_\infty,a}$. By Lemma~\ref{nodal smoothing} there exists
a $J_{n,t_\infty}$-holomorphic deformation of $\varphi'_{n,t_\infty}$
to an immersion $\varphi''_{n,t_\infty}$ containing $|\mathbf{x}^a|$
in the image. In the further deformation process we want to only
deform $\varphi'_{n,t_\infty}$ and keep the other components. To this
end change the path $J_{n,t}$ for $t>t_\infty$ so that it stays
constant on a neighbourhood of $\bigcup_{a'\neq a}
C_{n,t_\infty,a'}$. This is possible by Remark~\ref{generic J},2. We
also choose a new path $x'_n(t)$ so that $x'_n(t)\not\in
\bigcup_{a'\neq a} C_{n,t_\infty,a'}$. Now apply the previous
reasoning with $\varphi_n$ replaced by $\varphi''_{n,t_\infty}$ and
$\mathbf{x}_n(t)$ by the concatenation of $\mathbf{x}$ and $x'_n(t)$.

Because the degree $c_1(M)\cdot C_{n,t_\infty,a}$ decreases by an
integral amount each time we split the curve, after finitely many
changes of $\{J_{n,t}\}$ the deformation will be successful for all
$t$. The result is a family of cycles $C_{n,t}$ that is
pseudo-holomorphic for a family of $\mathscr{C}^{0,\alpha}$-small
perturbations of $J_n$ and containing a point $x_n(t)$ with
$$
d(x_n(t), |C_\infty|) \ge2\eps t.
$$
Finally, the inclusion $|\mathbf{x}|\subset |C_{n,t}|$ for all $t$
implies that $(m(C_{n,t}),\delta(C_{n,t}))$ is always smaller than
$(m(C_\infty), \delta(C_\infty))  =(m(C^0_{n,t}), \delta(C_{n,t})
+\delta^{\mathrm{par}})$. Note that a potential drop of $\delta$
happens on $B_\delta(|C^0_\infty|_\sing)$. Hence upon replacement of
the parasitic part $\delta$ decreases by $\delta^{\mathrm{par}}$, so
the induction hypothesis applies to the reglued curves. As above we
thus deduce that for any $t$ there exists a partial smoothing of
$C_{n,t}$ that is the image of a $J_{n,t}$-holomorphic map
$\varphi_{n,t}$ homotopic to $\varphi_n$ through a family of
pseudo-holomorphic immersions as demanded. Hence there exists $t\le
1$ with $\varphi'_n :=\varphi_{n,t}$ fulfilling the desired
properties~1--3.
\qed
\medskip

\noindent
7)\ \emph{The reduced case: Pseudo-holomorphic deformation of $C_n$
over a path $\{J_{n,t}\}_{t\in [0,1]}$ connecting $J_n$ with $J$.}\\
If $C_\infty$ is reduced there is no room for imposing one more point
constraint without spoiling unobstructedness of the deformation.
Instead choose general paths $\{J_{n,t}\}_{t\in [0,1]}$ of tamed almost
complex structures connecting $J_n$ with $J$ and such that
$$
J_{n,t}\stackrel{n\to\infty}{\lra} J
$$
uniformly in the $\mathscr{C}^0$-topology everywhere and in the
$\mathscr{C}^{0,\alpha}$-topology on $M\setminus A$. By
Remark~\ref{generic J},2 we can also arrange $J_{n,t}\in
\mathscr{J}_{|C_\infty|\cap U,\hol}$ for every $t>1/2$.

For $n\gg 0$ and $t\in[0,1)$ we seek to find a $J_{n,t}$-holomorphic
immersion $\varphi_{n,t}: (\Sigma_n,j_{n,t}) \to M$ with
\begin{itemize}
\item[(i)] $d_{\text{cyc}}(\varphi_{n,t}(\Sigma_n),C_\infty)<\eps$. 
\item[(ii)] $\varphi_{n,t}$ is homotopic to $\varphi_n$ through a family of
immersions, pseudo-holo\-mor\-phic with respect to a path of almost
complex structures staying $\eps$-close to $J_n$ in
$\mathscr{C}^{0}$ and in $\mathscr{C}^{0,\alpha}$ on $M\setminus
B_\eps(A)$.
\item[(iii)] $|\mathbf{x}|\subset \varphi_{n,t}(\Sigma_n)$.
\end{itemize}

\begin{lemma}\label{dichotomy}
For every $\eps>0$ one of the following two cases occurs.
\begin{enumerate}
\item There exists an $n>0$ so that
$\varphi_{n,t}$ with properties (i)--(iii) exists for all $t<1$.
\item For every $n\gg 0$ there exists $\varphi_{n,t_n}$ with properties
(i)--(iii) and so that
$$
d_{\text{cyc}}(\varphi_{n,t_n}(\Sigma_n),C_\infty)=\eps.
$$
\end{enumerate}
\end{lemma}

\proof
For $t=0$ we may take $\varphi_{n,t}=\varphi_n$. Let $\tau_n$ be
the maximal number so that for each $t\in[0,\tau_n)$ a map
$\varphi_{n,t}$ obeying (i)--(iii) exists. Assume that $\tau_n<1$. Let
$\{\varphi_{n,t_i}\}_i$ be a sequence of maps $\Sigma_n\to M$
obeying (i)--(iii) and with $t_i\nearrow \tau_n$ for $i\to\infty$.
By the Gromov compactness theorem, after going over to a subsequence,
we may assume cycle-theoretic convergence
$$
\varphi_{n,t_i}(\Sigma_n)\stackrel{i\to\infty}{\lra}
C_{\tau_n}.
$$
Since $\{J_{n,t}\}$ is generic at $\tau_n<1$ and
$|\mathbf{x}|\subset C_{\tau_n}$, Proposition~\ref{no equisingular
deformations of C_infty} implies
$$
(m(C_{\tau_n}),\delta(C_{\tau_n}))< (m(C_\infty),\delta(C_\infty)),
$$
unless $C_\infty$ is nodal. If $(m,\delta)$ drops argue as in the
non-reduced case: Regluing of the non-parasitic part and application
of the induction hypothesis shows the existence of $\varphi_{n,t}$
for $t>\tau_n$ with the requested properties, provided
$d_{\text{cyc}}(C_{\tau_n},C_\infty)<\eps$. Otherwise there exists
$\varphi_{n,t_n}$ as in (2) (with slightly smaller $\eps$). In the
nodal case the existence of $\varphi_{n,t}$ for $t>\tau_n$ follows
from the deformation theory of nodal curves Theorem~\ref{smoothing
nodal curves}.
\qed
\smallskip

Note that this line of reasoning fails in the non-reduced case,
because in the proof of closedness $|\mathbf{x}|$ may end up unevenly
distributed on $C_{n,\tau_n}$ if $|C_{n,\tau_n}|$ is reducible
near a multiple component of $C_\infty$. A smoothing of
$C_{n,\tau_n}$ preserving the incidence conditions may then not
exist.
\medskip

\noindent
8)\ \emph{Taking the limit.}\\
Assume first that either $C_\infty$ is non-reduced or that
Lemma~\ref{dichotomy},2 applies. Let $C'_n$ be the sequence of
deformations of $C_n$ constructed in this lemma or in
Lemma~\ref{C'_n} in Step~6 respectively. By the Gromov compactness
theorem we may assume convergence
$$
C'_n\stackrel{n\to\infty}{\lra} C'_\infty.
$$
By construction $d_\cyc (C'_\infty, C_\infty) \ge\eps$, hence
$C'_\infty\neq C_\infty$. Because also $|\mathbf{x}|\subset
C'_\infty$, Proposition~\ref{no equisingular deformations of C_infty}
implies 
$$
(m(C'_\infty),\delta (C'_\infty)) < (m(C_\infty),\delta(C_\infty)).
$$
At the expense of going over to a deformation ${J^0}'$ of
$J^0$ as allowed by Proposition~\ref{unobstructedness of cycle space}
replace the parasitic parts by the removed part of $C^0_\infty$, cf.\
Step~4. Take a ${J^0}'$-holomorphic smoothing $(C')^\dagger$ of this
curve, existent by Remark~\ref{smoothing exists}. The induction
hypothesis shows that $(C')^\dagger$ is symplectically isotopic to
$C'_n$ with the parasitic part replaced, hence to $C_n^0$.

In the remaining case Lemma~\ref{dichotomy},1 there is a family of
reglued curves $C^0_{n,t}$ that is pseudo-holomorphic for a family of
almost complex structures that are integrable in a fixed
neighbourhood of $|C_\infty|$ and fiberwise agreeing with $J$ for
$t>1/2$. In this case the isotopy statement follows for $n$
sufficiently large and $t$ close to $1$ from the unobstructedness
result Proposition~\ref{unobstructedness of cycle space}.

In the notation of the statement of the lemma this produces the
desired isotopy between $C_n$ and a smoothing of $C_\infty$ under the
additional hypothesis~($**$). This was left to be shown in Step~1.
The proof of Lemma~\ref{isotopy lemma} is finished.
\qed

%===========================================================
\section{Proofs of Theorems A, B and C}\noindent
This section provides the proofs of the three main theorems stated in
the introduction. We shall use one more Lemma.

\begin{lemma}\label{degree lemma}
Let $(M,J)$ be an almost complex manifold diffeomorphic to $\CC\PP^2$
or to an $S^2$-bundle over $S^2$. Assume that $C\subset M$ is a
non-trivial irreducible pseudo-holomorphic curve with $c_1(M)\cdot
C>0$. Then for all $m>1$
\begin{eqnarray}\label{degree estimate}
m(c_1(M)\cdot C-1)> c_1(M)\cdot C+g(C)-1
\end{eqnarray}
in either of the following cases.
\begin{enumerate}
\item[(i)] $M=\CC\PP^2$ and $C\sim dH$ with $d\le 8$, $H\subset
\CC\PP^2$ a hyperplane.
\item[(ii)] $M=\CC\PP^1\times\CC\PP^1$ or $\mathbb{F}_1$ and $C\sim
dH+kF$ with $d\le 3$, $k\ge0$, where $F$ and $H$ are a fiber and a
section of the symplectic $S^2$-bundle $M\to \CC\PP^1$ respectively
with $H\cdot H\in\{0,1\}$.
\item[(iii)] $M=\mathbb{F}_1$ and $C\sim dH+kF$ with $d+k\le 8$,
notations as in (ii).
\end{enumerate}
\end{lemma}

\proof
Since $c_1(M)\cdot C>0$ it suffices to establish
(\ref{degree estimate}) for $m=2$. For (i) the genus
formula~(\ref{genus formula}) yields
$$
g-1\le\frac{d^2-3d}{2}.
$$
Inequality (\ref{degree estimate}) thus follows from
$$
2(3d-1)>3d+\frac{d^2-3d}{2}\quad\text{or}\quad d^2-9d+4<0.
$$
This is true for $1\le d\le 8$.

For $M=\CC\PP^1\times\CC\PP^1$ it holds $H\cdot H=0$ and the genus
formula reads
$$
g-1\le \frac{(dH+kF)^2-(dH+kF)(2H+2F)}{2}
= dk-d-k.
$$
So (\ref{degree estimate}) will follow from
$$
2(2d+2k-1)> 2d+2k+(dk-d-k)=dk+d+k,
$$
which is equivalent to
$$
(k-3)(3-d)>-7.
$$
This is true for any $k> 0$ as long as $d\le 3$ and for $k=0$,
$1\le d\le 3$.

Finally the case $M=\mathbb{F}_1$. Then $H\cdot H=1$, $c_1(M)=2H+F$,
$$
g-1\le \frac{(dH+kF)^2-(dH+kF)(2H+F)}{2}
=\frac{d^2+2dk-3d-2k}{2},
$$
and (\ref{degree estimate}) follows from
$$
2(3d+2k-1)>3d+2k+\frac{d^2+2dk-3d-2k}{2}
=\frac{d^2+2dk+3d+2k}{2},
$$
or
\begin{eqnarray}\label{inequality F_1}
k(2d-6)+d^2-9d+4<0.
\end{eqnarray}
For $k\ge0$ and $1\le d\le 3$ the first term is non-positive and
$$
d^2-9d+4<0,
$$
so the result follows. For $d=0$, $k\ge 1$ we obtain
$$
-6k+4\le -2<0.
$$
If $d\ge 4$ but $k\le 8-d$ inequality (\ref{inequality F_1}) implies
as sufficient condition
$$
(8-d)(2d-6)+d^2-9d+4 = -d^2+13d-44 <0.
$$
This inequality holds true for all $d$. The proof is finished.
\qed

\begin{remark}
We gave (i) in the lemma only to illustrate the origin of the degree
restriction $d\le 17$ for $\CC\PP^2$. For technical
reasons we have to do the computation on the $\CC\PP^1$-bundle
$\mathbb{F}_1$, with some additional properties covered by
case (iii) in the lemma.
\end{remark}
\bigskip

\noindent
\emph{Proof of Theorem~B.}\ Let $\omega$ be the symplectic form on
$M$. Denote by $\Sigma\subset M$ the symplectic submanifold.  By
Proposition~\ref{fibered J} there exists a map $p:M\to\CC\PP^1$ and
an $\omega$-tamed almost complex structure $J$ on $M$ making $p$ a
pseudo-holomorphic map and $\Sigma$ a $J$-holomorphic curve. Let
$H,S,F$ be two disjoint sections and a fiber of $p$ respectively,
with $H\cdot H\in\{0,1\}$. According to Remark~\ref{wlog Sigma>0} we
may assume $c_1(M)\cdot \Sigma>0$. Then deformations of $\Sigma$ as
pseudo-holomorphic curve are unobstructed by the smooth case of
Theorem~\ref{smoothing nodal curves}.

Possibly after deforming $\Sigma$ slightly we may therefore assume
$J\in \mathscr{J}_\reg$ for the Baire subset $\mathscr{J}_\reg\subset
\mathscr{J}$ introduced in Section~\ref{generic paths}. Since
$\omega$ tames also the integrable complex structure $I$ on $M$,
Lemma~\ref{genericity c_1}, Theorem~\ref{genericity Seva} and
Proposition~\ref{very generic J} imply the existence of a path
$\{J_t\}_{t\in [0,1]}$ of $\omega$-tamed almost complex structures
having the following properties.
\begin{enumerate}
\item[(i)] $J_0=J$, $J_1=I$ and $p:(M,\omega,J_t)\to\CC\PP^1$ is
a symplectic pseudo-holomorphic $S^2$-bundle for every $t$.
\item[(ii)] $(M,J_t)$ is monotone for every $t$.
\item[(iii)] The path $\{J_t\}_{t\in [0,1]}$ is generic in the sense of
Theorem~\ref{genericity Seva},2.
\item[(iv)] Any $J_t$-holomorphic curve
$C\subset M$, $C\not\sim S$ is homologous to $dH+kF$ with $d\ge0$,
$k\ge 0$.
\end{enumerate}

Consider the set $\Omega$ of $\tau\in[0,1]$ so that for every $t\le\tau$ a
smooth $J_t$-holomorphic curve $C_t\subset M$ exists that is isotopic
to $\Sigma$ through an isotopy of $\omega$-symplectic submanifolds.
By definition $\Omega$ is an interval. It is open as subset of $[0,1]$ by
the smooth case of Theorem~\ref{smoothing nodal curves}. It remains
to show it is closed. Let $t_n\nearrow \tau$, $t_n\in \Omega$. For each
$n$ choose a $J_{t_n}$-holomorphic curve $C_n$ symplectically
isotopic to $\Sigma$. By the Gromov compactness theorem we may assume
that the $C_n$ converge to a $J_\tau$-holomorphic cycle $C_\infty
=\sum_a m_a C_{\infty,a}$.

Now $J_{t_n}\to J_\tau$, $C_n\to C_\infty$ fulfill the hypotheses of
Lemma~\ref{isotopy lemma}. Condition~($*$) follows from
Lemma~\ref{degree lemma},ii. The conclusion of Lemma~\ref{isotopy
lemma}  in conjunction with Remark~\ref{smoothing exists} now shows
that $\tau\in \Omega$. Hence $\Omega=[0,1]$.
(End of proof of Theorem~B.)
\qed
\medskip

\noindent
\emph{Proof of Theorem~C.}\
The proof is essentially the same as that of Theorem~B. As we are not
in a fibered situation we proceed as follows. Let $\sigma:
\mathbb{F}_1\to \CC\PP^2$ be the blowing up in a point $P\in
\CC\PP^2\setminus \Sigma$. Choose a K\"ahler form $\tilde\omega$ on
$\mathbb{F}_1$ that agrees with $\sigma^*(\omega)$ away from
$\sigma^{-1}(B_\eps(P))$ for $\eps< d(P,\Sigma)/2$. Now run the
program in the proof of Theorem~B to the $\tilde\omega$-symplectic
surface $\sigma^{-1}(\Sigma)$. By Remark~\ref{generic J} we are free
to take the $J_t$ and all intermediately occurring almost complex
structures to agree with the standard integrable complex structure on
$\sigma^{-1}(B_\eps(P))$. Then $J_t$ descends to an $\omega$-tamed
almost complex structure $\bar J_t$ on $\CC\PP^2$. If $C_t$ is a
smooth $J_t$-holomorphic curve homologous to $\sigma^{-1}(\Sigma)$
then $C_t$ is disjoint from $S=\sigma^{-1}(P)$ for homological
reasons. Thus $\sigma(C_t)$ is a smooth $\bar J_t$-holomorphic curve,
hence symplectic with respect to $\omega$.

It remains to verify assumption~($*$) in Lemma~\ref{isotopy lemma}.
Let $C'=\sum_{a\ge 0} m'_a C'_a$ be a $J_t$-holomorphic cycle
homologous to $dH$, where $H$ is a section with $H\cdot H=1$. By
Lemma~\ref{genericity c_1},1 there are no $J_t$-holomorphic curves
homologous to $aH+bF$ with $b<0$ except $S\sim H-F$. Therefore, if
$C'_a\sim d_a H+ k_a F$ with $k_a\ge 0$ for $a>0$, and $m'_0 C'_0= mS$
then
$$
d= m+\sum_{a>0} m'_a d_a\quadand m=\sum_{a>0} m'_ak_a.
$$
Hence $d=\sum_{a>0} m'_a(d_a+k_a)$. Now if $d\le 17$ and $C'$ has a
multiple component then $d_a+k_a\le 8$ for all $a$. Assumption~($*$)
therefore follows from Lemma~\ref{degree lemma},iii. (End of proof of
Theorem~C.)
\qed
\bigskip

\noindent
\emph{Proof of Theorem~A.}\ The main theorem of \cite{st99} implies
that $M\to S^2$ factors over a degree 2 cover of an $S^2$-bundle
$p:P\to S^2$, branched along a symplectic surface $B\subset P$ with
$p|_B$ a simply branched cover of $S^2$. The assumption on the
monodromy means that $B$ is connected. By Theorem~B there exists a
symplectic isotopy $B_t$ connecting $B=B_0$ to a holomorphic curve
$B_1$ with respect to the integrable complex structure. Apply
Proposition~\ref{fibered J} with $\Sigma_t=B_t$. The result is a
family $J_t$ of almost complex structures and a family of
$S^2$-bundles $p_t: P\to S^2$ with
\begin{enumerate}
\item $B_t$ is $J_t$-holomorphic.
\item $p_t:(P,J_t)\to \CC\PP^1$ is pseudo-holomorphic.
\end{enumerate}
Perturbing $B_t$ slightly using the smooth case of
Theorem~\ref{smoothing nodal curves} we can achieve that $B_t\to
\CC\PP^1$ has only simple branch points for every $t$. Let $M_t$ be
the two-fold cover of $P$ branched along $B_t$. Then $M_t\to S^2$ is
a genus-$2$ symplectic Lefschetz fibration for every $t$. This shows
that $M\to S^2$ is even isotopic to a holomorphic Lefschetz
fibration. (End of proof of Theorem~A.)
\qed

%===========================================================

\end{document}